# Quaternion Rhapsody

Aleks Kleyn

*Dedicated to I. M. Gelfand*

ABSTRACT. In this paper I explore the set of quaternion algebras over field. Quaternion algebra $E(C, -1, -1)$ is isomorphic to tensor product of complex field $C$ and quaternion algebra $H = E(R, -1, -1)$. Considered the set of quaternion functions, which satisfy to equation similar to Cauchy-Riemann equation for a complex function.

## Contents

## 1. Preface

When I started my research in area of calculus over division ring, I paid attention to large volume of papers dedicated to regular functions of quaternion and analogue of the Cauchy-Riemann equation. I wanted to understand whether the Cauchy-Riemann equation appears in the frame of the theory that I explore. This was the reason why I considered division ring as vector space over center.

Introduction of basis simplifies some constructions and presents a bridge between the Gâteaux derivative and Jacobian matrix of map. This is exactly the place, where the Cauchy-Riemann equations should appear. The exploration of derivative of function of complex numbers reveals that the Cauchy-Riemann equation has algebraic origin and is related with statement that there exists $R$-linear function over complex field, however this function is not $C$-linear. For instance, conjugation of complex number is linear over real field, however it is not linear map over complex field.

In quaternion algebra $H$, there exist only $R$-linear map. The corollary of this statement is the ability represent every $R$-linear map using quaternion and absence of the evident analogue of the Cauchy-Riemann equation in quaternion algebra. However, in numerous papers and books dedicated to calculus over quaternion algebra mathematicians explore sets of maps that have properties similar to properties of functions of complex variable. Some of authors do not restrict themselves to quaternion algebra and explore more general algebras.

In the paper [3], Gelfand explores the quaternion algebra over arbitrary field assuming that product depends on arbitrary parameters. We assume $H = E(R, -1, -1)$. Using this paper, I decided to explore two cases that are important for me.

The algebra $E(R, a, b)$ was interesting for me because I supposed to find parameters $a$, $b$ such that the system of linear equations [2]-(3.1.17) is singular. It was important to understand what happens in this case. When I explored the structure of linear map over division ring it was not evident how non singularity of system of linear equations [2]-(3.1.17) have an influence on answer. However, the solution to this problem was not the one I had expected. It turned out that the system of linear equations is so simple that everybody can see that this system cannot be singular.

I have wrote that the Cauchy-Riemann equation is related with statement that complex field has real field as subfield. I assumed also that similar statement is possible in algebras with enough aggregate center. This is why I expected to see analogue of the Cauchy-Riemann equation in the quaternion algebra over complex field. In the course of solving the problem I realized that algebra $E(C, -1, -1)$ is isomorphic to tensor product $C \otimes H$. Therefore linear functions of this algebra satisfy to the Cauchy-Riemann equation for $C$-component of tensor product. Therefore I can tell the same about Jacobian matrix of arbitrary function. Natural extension of this topic is exploration of tensor product $C \otimes (C \otimes H)$. In particular, I put my attention is non associativity of tensor product.

In the book [2] on the base of which I wrote this paper, I explore the Gâteaux derivative of function over division ring. However in this paper I consider arbitrary algebras that not always are division rings. During the time that I explore the Gâteaux derivative, I realized that this subject can be generalized to more wide set of algebras. I will prepare the complete research later. I wrote this paper in order to explore conditions when the Cauchy-Riemann equations are possible.

## 2. Conventions

(1) Function and mapping are synonyms. However according to tradition, correspondence between either rings or vector spaces is called mapping and a mapping of either real field or quaternion algebra is called function.
(2) We can consider division ring $D$ as $D$-vector space of dimension 1. According to this statement, we can explore not only homomorphisms of division ring $D_1$ into division ring $D_2$, but also linear maps of division rings. This means that map is multiplicative over maximum possible field. In particular, linear map of division ring $D$ is multiplicative over center $Z(D)$. This statement does not contradict with definition of linear map of field because for field $F$ is true $Z(F) = F$. When field $F$ is different from maximum possible, I explicit tell about this in text.

(3) Let $A$ be free finite dimensional algebra. Considering expansion of element of algebra $A$ relative basis $\overline{\overline{e}}$ we use the same root letter to denote this element and its coordinates. However we do not use vector notation in algebra. In expression $a^2$, it is not clear whether this is component of expansion of element $a$ relative basis, or this is operation $a^2 = aa$. To make text clearer we use separate color for index of element of algebra. For instance,
$$a = a^i \overline{e}_i$$

(4) If free finite dimensional algebra has unit, then we identify the vector of basis $\overline{e}_{\mathbf{0}}$ with unit of algebra.

(5) Without a doubt, the reader may have questions, comments, objections. I will appreciate any response.

## 3. Linear Function of Complex Field

**Theorem 3.1** (the Cauchy-Riemann equations). *Let us consider complex field $C$ as two-dimensional algebra over real field. Let*

(3.1) $$\overline{e}_{C \cdot \mathbf{0}} = 1 \quad \overline{e}_{C \cdot \mathbf{1}} = i$$

*be the basis of algebra $C$. Then in this basis product has form*

(3.2) $$\overline{e}_{C \cdot \mathbf{1}}^2 = -\overline{e}_{C \cdot \mathbf{0}}$$

*and structural constants have form*

(3.3) $$\begin{aligned} C_{C \cdot \mathbf{00}}^{\mathbf{0}} &= 1 & C_{C \cdot \mathbf{01}}^{\mathbf{1}} &= 1 \\ C_{C \cdot \mathbf{10}}^{\mathbf{1}} &= 1 & C_{C \cdot \mathbf{11}}^{\mathbf{0}} &= -1 \end{aligned}$$

*Matrix of linear function*
$$y^i = x^j f_j^i$$
*of complex field over real field satisfies relationship*

(3.4) $$f_{\mathbf{0}}^{\mathbf{0}} = f_{\mathbf{1}}^{\mathbf{1}}$$

(3.5) $$f_{\mathbf{0}}^{\mathbf{1}} = -f_{\mathbf{1}}^{\mathbf{0}}$$

*Proof.* Equations (3.2) and (3.3) follow from equation $i^2 = -1$. Using equation [2]-(3.1.17) we get relationships

(3.6) $$f_{\mathbf{0}}^{\mathbf{0}} = f^{kr} C_{C \cdot \mathbf{k0}}^{\mathbf{p}} C_{C \cdot \mathbf{pr}}^{\mathbf{0}} = f^{\mathbf{0}r} C_{C \cdot \mathbf{00}}^{\mathbf{0}} C_{C \cdot \mathbf{0}r}^{\mathbf{0}} + f^{\mathbf{1}r} C_{C \cdot \mathbf{10}}^{\mathbf{1}} C_{C \cdot \mathbf{1}r}^{\mathbf{0}} = f^{\mathbf{00}} - f^{\mathbf{11}}$$

(3.7) $$f_{\mathbf{0}}^{\mathbf{1}} = f^{kr} C_{C \cdot \mathbf{k0}}^{\mathbf{p}} C_{C \cdot \mathbf{pr}}^{\mathbf{1}} = f^{\mathbf{0}r} C_{C \cdot \mathbf{00}}^{\mathbf{0}} C_{C \cdot \mathbf{0}r}^{\mathbf{1}} + f^{\mathbf{1}r} C_{C \cdot \mathbf{10}}^{\mathbf{1}} C_{C \cdot \mathbf{1}r}^{\mathbf{1}} = f^{\mathbf{01}} + f^{\mathbf{10}}$$

(3.8) $$f_{\mathbf{1}}^{\mathbf{0}} = f^{kr} C_{C \cdot \mathbf{k1}}^{\mathbf{p}} C_{C \cdot \mathbf{pr}}^{\mathbf{0}} = f^{\mathbf{0}r} C_{C \cdot \mathbf{01}}^{\mathbf{1}} C_{C \cdot \mathbf{1}r}^{\mathbf{0}} + f^{\mathbf{1}r} C_{C \cdot \mathbf{11}}^{\mathbf{0}} C_{C \cdot \mathbf{0}r}^{\mathbf{0}} = -f^{\mathbf{01}} - f^{\mathbf{10}}$$

(3.9) $$f_{\mathbf{1}}^{\mathbf{1}} = f^{kr} C_{C \cdot \mathbf{k1}}^{\mathbf{p}} C_{C \cdot \mathbf{pr}}^{\mathbf{1}} = f^{\mathbf{0}r} C_{C \cdot \mathbf{01}}^{\mathbf{1}} C_{C \cdot \mathbf{1}r}^{\mathbf{1}} + f^{\mathbf{1}r} C_{C \cdot \mathbf{11}}^{\mathbf{0}} C_{C \cdot \mathbf{0}r}^{\mathbf{1}} = f^{\mathbf{00}} - f^{\mathbf{11}}$$

(3.4) follows from equations (3.6) and (3.9). (3.5) follows from equations (3.7) and (3.8). □

**Theorem 3.2** (the Cauchy-Riemann equations)**.** *Since matrix*

$$\begin{pmatrix} \dfrac{\partial y^0}{\partial x^0} & \dfrac{\partial y^0}{\partial x^1} \\ \dfrac{\partial y^1}{\partial x^0} & \dfrac{\partial y^1}{\partial x^1} \end{pmatrix}$$

*is Jacobian matrix of map of complex variable*

$$x = x^0 + x^1 i \to y = y^0(x^0, x^1) + y^1(x^0, x^1)i$$

*over real field, then*

(3.10)
$$\begin{aligned} \dfrac{\partial y^1}{\partial x^0} &= -\dfrac{\partial y^0}{\partial x^1} \\ \dfrac{\partial y^0}{\partial x^0} &= \dfrac{\partial y^1}{\partial x^1} \end{aligned}$$

*Proof.* The statement of theorem is corollary of theorem 3.1. □

**Theorem 3.3.** *Derivative of function of complex variable satisfyes to equation*

(3.11)
$$\dfrac{\partial y}{\partial x^0} + i\dfrac{\partial y}{\partial x^1} = 0$$

*Proof.* Equation

$$\dfrac{\partial y^0}{\partial x^0} + i\dfrac{\partial y^1}{\partial x^0} + i\dfrac{\partial y^0}{\partial x^1} - \dfrac{\partial y^1}{\partial x^1} = 0$$

follows from equations (3.10). □

Equation (3.11) is equivalent to equation

(3.12)
$$\begin{pmatrix} 1 & i \end{pmatrix} \begin{pmatrix} \dfrac{\partial y^0}{\partial x^0} & \dfrac{\partial y^0}{\partial x^1} \\ \dfrac{\partial y^1}{\partial x^0} & \dfrac{\partial y^1}{\partial x^1} \end{pmatrix} \begin{pmatrix} 1 \\ i \end{pmatrix} = 0$$

4. Quaternion Algebra

In this paper I explore the set of quaternion algebras defined in [3].

**Definition 4.1.** Let $F$ be field. Extension field $F(i, j, k)$ is called **the quaternion algebra** $E(F, a, b)$ **over the field** $F$[1] if multiplication in algebra $E$ is defined according to rule

(4.1)

|   | $i$ | $j$ | $k$ |
|---|-----|-----|-----|
| $i$ | $a$ | $k$ | $aj$ |
| $j$ | $-k$ | $b$ | $-bi$ |
| $k$ | $-aj$ | $bi$ | $-ab$ |

where $a, b \in F$, $ab \neq 0$. □

---

[1] I follow definition from [3].

Elements of the algebra $E(F, a, b)$ have form

$$x = x^0 + x^1 i + x^2 j + x^3 k$$

where $x^i \in F$, $i = 0, 1, 2, 3$. Quaternion

$$\overline{x} = x^0 - x^1 i - x^2 j - x^3 k$$

is called conjugate to the quaternion $x$. We define **the norm of the quaternion** $x$ using equation

(4.2) $$|x|^2 = x\overline{x} = (x^0)^2 - a(x^1)^2 - b(x^2)^2 + ab(x^3)^2$$

From equation (4.2), it follows that $E(F, a, b)$ is algebra with division only when $a < 0$, $b < 0$. In this case we can renorm basis such that $a = -1$, $b = -1$.

We use symbol $E(F)$ to denote the quaternion division algebra $E(F, -1, -1)$ over the field $F$. We will use notation $H = E(R, -1, -1)$. Multiplication in quaternion algebra $E(F)$ is defined according to rule

(4.3)
$$\begin{array}{c|ccc} & i & j & k \\ \hline i & -1 & k & -j \\ j & -k & -1 & i \\ k & j & -i & -1 \end{array}$$

In algebra $E(F)$, the norm of the quaternion has form

(4.4) $$|x|^2 = x\overline{x} = (x^0)^2 + (x^1)^2 + (x^2)^2 + (x^3)^2$$

In this case inverse element has form

(4.5) $$x^{-1} = |x|^{-2}\overline{x}$$

The inner automorphism of quaternion algebra $H$[2]

(4.6)
$$p \to qpq^{-1}$$
$$q(ix + jy + kz)q^{-1} = ix' + jy' + kz'$$

describes the rotation of the vector with coordinates $x$, $y$, $z$. The norm of quaternion $q$ is irrelevant, although usually we assume $|q| = 1$. If $q$ is written as sum of scalar and vector

$$q = \cos\alpha + (ia + jb + kc)\sin\alpha \quad a^2 + b^2 + c^2 = 1$$

then (4.6) is a rotation of the vector $(x, y, z)$ about the vector $(a, b, c)$ through an angle $2\alpha$.

## 5. Tower of Algebras

Let $F_1$ be algebra over the field $F_2$. Let $\overline{\overline{e}}_{12}$ be basis of algebra $F_1$ over the field $F_2$. Let $C_{12 \cdot ij}^{\phantom{12\cdot}k}$ be structural constants of algebra $F_1$ over the field $F_2$.

Let $F_2$ be algebra over the field $F_3$. Let $\overline{\overline{e}}_{23}$ be basis of algebra $F_2$ over the field $F_3$. Let $C_{23 \cdot ij}^{\phantom{23\cdot}k}$ be structural constants of algebra $F_2$ over the field $F_3$.

I will consider the algebra $F_1$ as direct sum of algebras $F_2$. Each item of sum I identify with vector of basis $\overline{\overline{e}}_{12}$. Accordingly, I can consider algebra $F_1$ as algebra over field $F_3$. Let $\overline{\overline{e}}_{13}$ be basis of algebra $F_1$ over the field $F_3$. Index of basis $\overline{\overline{e}}_{13}$ consists from two indexes: index of fiber and index of vector of basis $\overline{\overline{e}}_{23}$ in fiber.

---

[2]See [6], p.643.

I will identify vector of basis $\overline{e}_{12 \cdot i}$ with unit in corresponding fiber. Then

(5.1) $$\overline{e}_{13 \cdot ji} = \overline{e}_{23 \cdot j} \overline{e}_{12 \cdot i}$$

The product of vectors of basis $\overline{\overline{e}}_{13}$ has form

(5.2) $$\overline{e}_{13 \cdot ji} \overline{e}_{13 \cdot mk} = \overline{e}_{23 \cdot j} \overline{e}_{12 \cdot i} \overline{e}_{23 \cdot m} \overline{e}_{12 \cdot k} = C_{23 \cdot jm}^{\ \ a} \overline{e}_{23 \cdot a} C_{12 \cdot ik}^{\ \ b} \overline{e}_{12 \cdot b}$$

Because $C_{12 \cdot ik}^{\ \ b} \in F_2$, then expansion $C_{12 \cdot ik}^{\ \ b}$ relative to basis $\overline{\overline{e}}_{23}$ has form

(5.3) $$C_{12 \cdot ik}^{\ \ b} = C_{12 \cdot ik}^{\ \ bc} \overline{e}_{23 \cdot c}$$

Let us substitute (5.3) into (5.2)

(5.4) $$\begin{aligned}\overline{e}_{13 \cdot ji} \overline{e}_{13 \cdot mk} &= C_{23 \cdot jm}^{\ \ a} \overline{e}_{23 \cdot a} C_{12 \cdot ik}^{\ \ bc} \overline{e}_{23 \cdot c} \overline{e}_{12 \cdot b} \\ &= C_{23 \cdot jm}^{\ \ a} C_{23 \cdot ac}^{\ \ d} C_{12 \cdot ik}^{\ \ bc} \overline{e}_{23 \cdot d} \overline{e}_{12 \cdot b} \\ &= C_{23 \cdot jm}^{\ \ a} C_{23 \cdot ac}^{\ \ d} C_{12 \cdot ik}^{\ \ bc} \overline{e}_{13 \cdot db}\end{aligned}$$

Therefore, we can define structural constants of algebra $F_1$ over field $F_3$

(5.5) $$C_{13 \cdot \cdot ji \cdot mk}^{\ \ \ \cdot db} = C_{23 \cdot jm}^{\ \ a} C_{23 \cdot ac}^{\ \ d} C_{12 \cdot ik}^{\ \ bc}$$

To verify construction, let us consider the product

(5.6) $$\begin{aligned}\overline{e}_{12 \cdot i} \overline{e}_{12 \cdot k} &= \overline{e}_{13 \cdot 0i} \overline{e}_{13 \cdot 0k} \\ &= C_{13 \cdot \cdot 0i \cdot 0k}^{\ \ \ \cdot db} \overline{e}_{13 \cdot db} \\ &= C_{23 \cdot 00}^{\ \ a} C_{23 \cdot ac}^{\ \ d} C_{12 \cdot ik}^{\ \ bc} \overline{e}_{13 \cdot db} \\ &= C_{23 \cdot 00}^{\ \ 0} C_{23 \cdot 0c}^{\ \ d} C_{12 \cdot ik}^{\ \ bc} \overline{e}_{13 \cdot bd} \\ &= C_{23 \cdot 0c}^{\ \ d} C_{12 \cdot ik}^{\ \ bc} \overline{e}_{13 \cdot bd}\end{aligned}$$

On the other hand

(5.7) $$\begin{aligned}\overline{e}_{12 \cdot i} \overline{e}_{12 \cdot k} &= C_{12 \cdot ik}^{\ \ b} \overline{e}_{12 \cdot b} \\ &= C_{12 \cdot ik}^{\ \ bd} \overline{e}_{23 \cdot d} \overline{e}_{23 \cdot 0} \overline{e}_{12 \cdot b} \\ &= C_{12 \cdot ik}^{\ \ bd} C_{23 \cdot 0d}^{\ \ c} \overline{e}_{23 \cdot c} \overline{e}_{12 \cdot b} \\ &= C_{12 \cdot ik}^{\ \ bd} C_{23 \cdot 0d}^{\ \ c} \overline{e}_{13 \cdot bc}\end{aligned}$$

Expressions (5.6), (5.7) coincide.

**Theorem 5.1.** *If $C_{12 \cdot ik}^{\ \ b} \in F_3$, then we multiply components $\overline{e}_{12 \cdot i}$, $\overline{e}_{23 \cdot k}$ of vector $\overline{e}_{13 \cdot ik}$ independently*

(5.8) $$C_{13 \cdot \cdot ji \cdot mk}^{\ \ \ \cdot db} = C_{23 \cdot jm}^{\ \ d} C_{12 \cdot ik}^{\ \ b}$$

*Proof.* From equation (5.2), it follows

(5.9) $$\overline{e}_{13 \cdot ji} \overline{e}_{13 \cdot mk} = C_{23 \cdot jm}^{\ \ a} \overline{e}_{23 \cdot a} C_{12 \cdot ik}^{\ \ b} \overline{e}_{12 \cdot b} = C_{23 \cdot jm}^{\ \ a} C_{12 \cdot ik}^{\ \ b} \overline{e}_{23 \cdot ab}$$

Equation (5.8) follows from the equation (5.9). $\square$

**Theorem 5.2.** *Let $F_1$ be algebra over the field $F_2$. Let $F_2$ be algebra over the field $F_3$. There exists equation*

(5.10) $$a_{12}^k = a_{13}^{ik} \overline{e}_{23 \cdot i}$$

*between coordinates $a_{12}^k$ of element $a \in F_1$ relative to basis $\overline{\overline{e}}_{12}$ and coordinates $a_{13}^{jk}$ of element $a \in F_1$ relative to basis $\overline{\overline{e}}_{13}$*

*Proof.* Element $a \in F_1$ has expansion

(5.11) $$a = a_{12}^{k} \overline{\overline{e}}_{12 \cdot k}$$

relative basis $\overline{\overline{e}}_{12}$ and expansion

(5.12) $$a = a_{13}^{ik} \overline{\overline{e}}_{13 \cdot ik}$$

relative basis $\overline{\overline{e}}_{13}$ From equations (5.1), (5.12), it follows that

(5.13) $$a = a_{13}^{ik} \overline{\overline{e}}_{23 \cdot i} \overline{\overline{e}}_{12 \cdot k}$$

The equation (5.10) follows from equations (5.11), (5.13). □

**Theorem 5.3** (The Cauchy-Riemann equations). *Let $F_1$ be algebra over the field $F_2$. Let $F_2$ be algebra over the field $F_3$. Matrix of $F_2$-linear mapping*

$$f : F_1 \to F_1$$

*of $F_2$-algebra $F_1$ satisfies relationship*

(5.14) $$f_{jl}^{ik} C_{23 \cdot mp}^{\phantom{23 \cdot}j} = C_{23 \cdot mn}^{\phantom{23 \cdot}i} f_{pl}^{nk}$$

*Proof.* Since $f$ is $F_2$-linear mapping, then

(5.15) $$f(ax) = af(x) \quad a \in F_2 \quad x \in F_1$$

Let

$$a = a^j \overline{e}_{23 \cdot j} \quad x = x_{13}^{jl} \overline{\overline{e}}_{13 \cdot jl}$$

Then

(5.16) $$f(ax)^{ik} = (ax)^{jl} f_{jl}^{ik} = C_{23 \cdot mp}^{\phantom{23 \cdot}j} a^m x^{pl} f_{jl}^{ik}$$

Similarly,

(5.17) $$(af(x))^{ik} = C_{23 \cdot mn}^{\phantom{23 \cdot}i} a^m (f(x))^{nk} = C_{23 \cdot mn}^{\phantom{23 \cdot}i} a^m f_{pl}^{nk} x^{pl}$$

From equations (5.15), (5.16), (5.17), it follows that

(5.18) $$f_{jl}^{ik} C_{23 \cdot mp}^{\phantom{23 \cdot}j} a^m x^{pl} = C_{23 \cdot mn}^{\phantom{23 \cdot}i} a^m f_{pl}^{nk} x^{pl}$$

Since $a \in F_2$, $x \in F_1$ are arbitrary, then the equation (5.14) follows from equation (5.18) □

**Theorem 5.4** (The Cauchy-Riemann equations). *Matrix of $F_2$-linear mapping*

$$f : F_2 \to F_2$$

*of $F_3$-algebra $F_2$ satisfies relationship*

(5.19) $$f_j^i C_{23 \cdot mp}^{\phantom{23 \cdot}j} = C_{23 \cdot mn}^{\phantom{23 \cdot}i} f_p^n$$

*Proof.* Statement of the theorem is proved if we assume $F_1 = F_2$ in the theorem 5.3. □

**Theorem 5.5** (The Cauchy-Riemann equations). *Matrix of $C$-linear mapping*

$$f : C \to C$$

*of $R$-algebra $C$ satisfies relationship*

(5.20) $$\begin{aligned} f_0^0 &= f_1^1 \\ f_1^0 &= -f_0^1 \end{aligned}$$

*Proof.* For complex field $C$, the equation (5.19) has form

$$(5.21) \quad \begin{array}{l} f_0^0 C_{C \cdot mp}^{\phantom{C}0} + f_1^0 C_{C \cdot mp}^{\phantom{C}1} = C_{C \cdot m0}^{\phantom{C}0} f_p^0 + C_{C \cdot m1}^{\phantom{C}0} f_p^1 \\ f_0^1 C_{C \cdot mp}^{\phantom{C}0} + f_1^1 C_{C \cdot mp}^{\phantom{C}1} = C_{C \cdot m0}^{\phantom{C}1} f_p^0 + C_{C \cdot m1}^{\phantom{C}1} f_p^1 \end{array}$$

From equations (3.3), (5.21), it follows that

$$(5.22) \quad \begin{cases} f_0^0 C_{C \cdot 00}^{\phantom{C}0} = C_{C \cdot 00}^{\phantom{C}0} f_0^0 \\ f_0^0 C_{C \cdot 11}^{\phantom{C}0} = C_{C \cdot 11}^{\phantom{C}0} f_1^1 \\ f_1^0 C_{C \cdot 01}^{\phantom{C}1} = C_{C \cdot 00}^{\phantom{C}0} f_1^0 \\ f_1^0 C_{C \cdot 10}^{\phantom{C}1} = C_{C \cdot 11}^{\phantom{C}0} f_0^1 \\ f_0^1 C_{C \cdot 00}^{\phantom{C}0} = C_{C \cdot 01}^{\phantom{C}1} f_0^1 \\ f_0^1 C_{C \cdot 11}^{\phantom{C}0} = C_{C \cdot 10}^{\phantom{C}1} f_1^0 \\ f_1^1 C_{C \cdot 01}^{\phantom{C}1} = C_{C \cdot 01}^{\phantom{C}1} f_1^1 \\ f_1^1 C_{C \cdot 10}^{\phantom{C}1} = C_{C \cdot 10}^{\phantom{C}1} f_0^0 \end{cases}$$

Equations (5.20) follow from equations (3.3), (5.22). □

I intentionally wrote the equation (5.22) in the intermediate form, to show how much more is tedious proof of the theorem 5.5 compared with proof of the analogous theorem 3.1. never used the statement that the product in the $F_1$-algebra $F_2$ is commutative. Therefore, the theorem 5.4 can be applied to an arbitrary algebra.

**Theorem 5.6** (The Cauchy-Riemann equations). *Matrix of $H$-linear mapping*

$$f : H \to H$$

*of $R$-algebra $H$ satisfies relationship*

$$(5.23) \quad \begin{array}{rrrr} f_0^0 = & f_1^1 = & f_2^2 = & f_3^3 \\ f_1^0 = -f_0^1 = -f_3^2 = & f_2^3 \\ f_2^0 = & f_3^1 = -f_0^2 = -f_1^3 \\ f_3^0 = -f_2^1 = & f_1^2 = -f_0^3 \end{array}$$

*Proof.* According to the theorem [2]-7.3.1, structural constants of quaternion algebra have form

$$(5.24) \quad \begin{array}{llll} C_{H \cdot 00}^{\phantom{H}0} = 1 & C_{H \cdot 01}^{\phantom{H}1} = \phantom{-}1 & C_{H \cdot 02}^{\phantom{H}2} = \phantom{-}1 & C_{H \cdot 03}^{\phantom{H}3} = \phantom{-}1 \\ C_{H \cdot 10}^{\phantom{H}1} = 1 & C_{H \cdot 11}^{\phantom{H}0} = -1 & C_{H \cdot 12}^{\phantom{H}3} = \phantom{-}1 & C_{H \cdot 13}^{\phantom{H}2} = -1 \\ C_{H \cdot 20}^{\phantom{H}2} = 1 & C_{H \cdot 21}^{\phantom{H}3} = -1 & C_{H \cdot 22}^{\phantom{H}0} = -1 & C_{H \cdot 23}^{\phantom{H}1} = \phantom{-}1 \\ C_{H \cdot 30}^{\phantom{H}3} = 1 & C_{H \cdot 31}^{\phantom{H}2} = \phantom{-}1 & C_{H \cdot 32}^{\phantom{H}1} = -1 & C_{H \cdot 33}^{\phantom{H}0} = -1 \end{array}$$

For quaternion algebra $H$, the equation (5.19) has form

$$(5.25) \quad \begin{array}{l} f_0^0 C_{H \cdot mp}^{\phantom{H}0} + f_1^0 C_{H \cdot mp}^{\phantom{H}1} + f_2^0 C_{H \cdot mp}^{\phantom{H}2} + f_3^0 C_{H \cdot mp}^{\phantom{H}3} \\ = C_{H \cdot m0}^{\phantom{H}0} f_p^0 + C_{H \cdot m1}^{\phantom{H}0} f_p^1 + C_{H \cdot m2}^{\phantom{H}0} f_p^2 + C_{H \cdot m3}^{\phantom{H}0} f_p^3 \quad i = 0 \\ f_0^1 C_{H \cdot mp}^{\phantom{H}0} + f_1^1 C_{H \cdot mp}^{\phantom{H}1} + f_2^1 C_{H \cdot mp}^{\phantom{H}2} + f_3^1 C_{H \cdot mp}^{\phantom{H}3} \\ = C_{H \cdot m0}^{\phantom{H}1} f_p^0 + C_{H \cdot m1}^{\phantom{H}1} f_p^1 + C_{H \cdot m2}^{\phantom{H}1} f_p^2 + C_{H \cdot m3}^{\phantom{H}1} f_p^3 \quad i = 1 \\ f_0^2 C_{H \cdot mp}^{\phantom{H}0} + f_1^2 C_{H \cdot mp}^{\phantom{H}1} + f_2^2 C_{H \cdot mp}^{\phantom{H}2} + f_3^2 C_{H \cdot mp}^{\phantom{H}3} \\ = C_{H \cdot m0}^{\phantom{H}2} f_p^0 + C_{H \cdot m1}^{\phantom{H}2} f_p^1 + C_{H \cdot m2}^{\phantom{H}2} f_p^2 + C_{H \cdot m3}^{\phantom{H}2} f_p^3 \quad i = 2 \\ f_0^3 C_{H \cdot mp}^{\phantom{H}0} + f_1^3 C_{H \cdot mp}^{\phantom{H}1} + f_2^3 C_{H \cdot mp}^{\phantom{H}2} + f_3^3 C_{H \cdot mp}^{\phantom{H}3} \\ = C_{H \cdot m0}^{\phantom{H}3} f_p^0 + C_{H \cdot m1}^{\phantom{H}3} f_p^1 + C_{H \cdot m2}^{\phantom{H}3} f_p^2 + C_{H \cdot m3}^{\phantom{H}3} f_p^3 \quad i = 3 \end{array}$$

From equations (5.24), (5.25), it follows that[3]

$$(5.26) \quad \begin{cases} f_1^0 C_{H \cdot 10}^{\phantom{H \cdot}1} = C_{H \cdot 11}^{\phantom{H \cdot}0} f_0^1 & i=0 \ \ m=1 \ \ p=0 \\ f_0^0 C_{H \cdot 11}^{\phantom{H \cdot}0} = C_{H \cdot 11}^{\phantom{H \cdot}0} f_1^1 & i=0 \ \ m=1 \ \ p=1 \\ f_3^0 C_{H \cdot 12}^{\phantom{H \cdot}3} = C_{H \cdot 11}^{\phantom{H \cdot}0} f_2^1 & i=0 \ \ m=1 \ \ p=2 \\ f_2^0 C_{H \cdot 13}^{\phantom{H \cdot}2} = C_{H \cdot 11}^{\phantom{H \cdot}0} f_3^1 & i=0 \ \ m=1 \ \ p=3 \\ f_2^0 C_{H \cdot 20}^{\phantom{H \cdot}2} = C_{H \cdot 22}^{\phantom{H \cdot}0} f_0^2 & i=0 \ \ m=2 \ \ p=0 \\ f_3^0 C_{H \cdot 21}^{\phantom{H \cdot}3} = C_{H \cdot 22}^{\phantom{H \cdot}0} f_1^2 & i=0 \ \ m=2 \ \ p=1 \\ f_0^0 C_{H \cdot 22}^{\phantom{H \cdot}0} = C_{H \cdot 22}^{\phantom{H \cdot}0} f_2^2 & i=0 \ \ m=2 \ \ p=2 \\ f_1^0 C_{H \cdot 23}^{\phantom{H \cdot}1} = C_{H \cdot 22}^{\phantom{H \cdot}0} f_3^2 & i=0 \ \ m=2 \ \ p=3 \\ f_3^0 C_{H \cdot 30}^{\phantom{H \cdot}3} = C_{H \cdot 33}^{\phantom{H \cdot}0} f_0^3 & i=0 \ \ m=3 \ \ p=0 \\ f_2^0 C_{H \cdot 31}^{\phantom{H \cdot}2} = C_{H \cdot 33}^{\phantom{H \cdot}0} f_1^3 & i=0 \ \ m=3 \ \ p=1 \\ f_1^0 C_{H \cdot 32}^{\phantom{H \cdot}1} = C_{H \cdot 33}^{\phantom{H \cdot}0} f_2^3 & i=0 \ \ m=3 \ \ p=2 \\ f_0^0 C_{H \cdot 33}^{\phantom{H \cdot}0} = C_{H \cdot 33}^{\phantom{H \cdot}0} f_3^3 & i=0 \ \ m=3 \ \ p=3 \end{cases}$$

Equations (5.23) follow from equations (5.24), (5.26). $\square$

**Theorem 5.7.** *Let $F_1$ be algebra over the field $F_2$. Let $F_2$ be algebra over the field $F_3$. Let $F_2$-linear mapping*

$$f : F_1 \to F_1$$

*have representation*

$$(5.27) \qquad f(a)^i = f_j^i a^j \quad f_j^i = f_j^{ik} \overline{\overline{e}}_{23 \cdot k}$$

*relative to basis $\overline{\overline{e}}_{12}$ and have representation*

$$(5.28) \qquad f(a)^{ki} = f_{lj}^{ki} a^{lj}$$

*relative to basis $\overline{\overline{e}}_{13}$. Coordinates $f_{il}^{jk}$ of mapping $f$ relative to basis $\overline{\overline{e}}_{13}$ have form*

$$(5.29) \qquad f_{lj}^{pi} = C_{23 \cdot kl}^{\phantom{23 \cdot}p} f_j^{ik}$$

*Proof.* Since

$$a^j = a^{lj} \overline{\overline{e}}_{23 \cdot l}$$

then, from the equation (5.27), it follows that

$$(5.30) \qquad f(a)^i = C_{23 \cdot kl}^{\phantom{23 \cdot}p} f_j^{ik} a^{lj} \overline{\overline{e}}_{23 \cdot p}$$

From the equation (5.28), it follows that

$$(5.31) \qquad f(a)^i = f(a)^{pi} \overline{\overline{e}}_{23 \cdot p} = f_{lj}^{pi} a^{lj} \overline{\overline{e}}_{23 \cdot p}$$

From equations (5.30), (5.31), it follows that

$$(5.32) \qquad f_{lj}^{pi} a^{lj} \overline{\overline{e}}_{23 \cdot p} = C_{23 \cdot kl}^{\phantom{23 \cdot}p} f_j^{ik} a^{lj} \overline{\overline{e}}_{23 \cdot p}$$

The equation (5.29) follows from the equation (5.32) because $\overline{\overline{e}}_{23}$ is basis and coordinates $a^{lj}$ of element $a \in F_1$ are arbitrary. $\square$

---

[3] I do not consider equations which are trivial or follow from the equations (5.26).

## 6. Quaternion Algebra over Complex Field

In this section, I will consider quaternion algebra $E(C,-1,-1)$, where $C$ is complex field.

Product in algebra $E(C,-1,-1)$ is defined according to table

(6.1)

| | $\overline{e}_{12\cdot\mathbf{1}}$ | $\overline{e}_{12\cdot\mathbf{2}}$ | $\overline{e}_{12\cdot\mathbf{3}}$ |
|---|---|---|---|
| $\overline{e}_{12\cdot\mathbf{1}}$ | $-1$ | $\overline{e}_{12\cdot\mathbf{3}}$ | $-\overline{e}_{12\cdot\mathbf{2}}$ |
| $\overline{e}_{12\cdot\mathbf{2}}$ | $-\overline{e}_{12\cdot\mathbf{3}}$ | $-1$ | $\overline{e}_{12\cdot\mathbf{1}}$ |
| $\overline{e}_{12\cdot\mathbf{3}}$ | $\overline{e}_{12\cdot\mathbf{2}}$ | $-\overline{e}_{12\cdot\mathbf{1}}$ | $-1$ |

According to theorem [2]-7.3.1, structural constants of quaternion algebra have form

$$C_{12\cdot\mathbf{00}}^{\mathbf{0}} = 1 \quad C_{12\cdot\mathbf{01}}^{\mathbf{1}} = 1 \quad C_{12\cdot\mathbf{02}}^{\mathbf{2}} = 1 \quad C_{12\cdot\mathbf{03}}^{\mathbf{3}} = 1$$
$$C_{12\cdot\mathbf{10}}^{\mathbf{1}} = 1 \quad C_{12\cdot\mathbf{11}}^{\mathbf{0}} = -1 \quad C_{12\cdot\mathbf{12}}^{\mathbf{3}} = 1 \quad C_{12\cdot\mathbf{13}}^{\mathbf{2}} = -1$$
$$C_{12\cdot\mathbf{20}}^{\mathbf{2}} = 1 \quad C_{12\cdot\mathbf{21}}^{\mathbf{3}} = -1 \quad C_{12\cdot\mathbf{22}}^{\mathbf{0}} = -1 \quad C_{12\cdot\mathbf{23}}^{\mathbf{1}} = 1$$
$$C_{12\cdot\mathbf{30}}^{\mathbf{3}} = 1 \quad C_{12\cdot\mathbf{31}}^{\mathbf{2}} = 1 \quad C_{12\cdot\mathbf{32}}^{\mathbf{1}} = -1 \quad C_{12\cdot\mathbf{33}}^{\mathbf{0}} = -1$$

Let $\overline{e}_{23\cdot k} = \overline{e}_{C\cdot k}$. Product in algebra $C$ is defined according to rule (3.2). According to the theorem 3.1, structural constants of complex field over real field have form (3.3)

Therefore, algebra $E(C,-1,-1)$ is isomorphic to tensor product $C \otimes H$. So we can select basis

$$\overline{e}_{13\cdot ij} = \overline{e}_{C\cdot i} \otimes \overline{e}_{12\cdot j}$$

(6.2)
$$\overline{e}_{13\cdot\cdot\mathbf{00}} = 1 \otimes 1 \quad \overline{e}_{13\cdot\cdot\mathbf{01}} = 1 \otimes i \quad \overline{e}_{13\cdot\cdot\mathbf{02}} = 1 \otimes j \quad \overline{e}_{13\cdot\cdot\mathbf{03}} = 1 \otimes k$$
$$\overline{e}_{13\cdot\cdot\mathbf{10}} = i \otimes 1 \quad \overline{e}_{13\cdot\cdot\mathbf{11}} = i \otimes i \quad \overline{e}_{13\cdot\cdot\mathbf{12}} = i \otimes j \quad \overline{e}_{13\cdot\cdot\mathbf{13}} = i \otimes k$$

**Theorem 6.1.** *Table of product in algebra $E(C,-1,-1)$ over field $R$ has form*

| | $1 \otimes i$ | $1 \otimes j$ | $1 \otimes k$ | $i \otimes 1$ | $i \otimes i$ | $i \otimes j$ | $i \otimes k$ |
|---|---|---|---|---|---|---|---|
| $1 \otimes i$ | $-1 \otimes 1$ | $1 \otimes k$ | $-1 \otimes j$ | $i \otimes i$ | $-i \otimes 1$ | $i \otimes k$ | $-i \otimes j$ |
| $1 \otimes j$ | $-1 \otimes k$ | $-1 \otimes 1$ | $1 \otimes i$ | $i \otimes j$ | $-i \otimes k$ | $-i \otimes 1$ | $i \otimes i$ |
| $1 \otimes k$ | $1 \otimes j$ | $-1 \otimes i$ | $-1 \otimes 1$ | $i \otimes k$ | $i \otimes j$ | $-i \otimes i$ | $-i \otimes 1$ |
| $i \otimes 1$ | $i \otimes i$ | $i \otimes j$ | $i \otimes k$ | $-1 \otimes 1$ | $-1 \otimes i$ | $-1 \otimes j$ | $-1 \otimes k$ |
| $i \otimes i$ | $-i \otimes 1$ | $i \otimes k$ | $-i \otimes j$ | $-1 \otimes i$ | $1 \otimes 1$ | $-1 \otimes k$ | $1 \otimes j$ |
| $i \otimes j$ | $-i \otimes k$ | $-i \otimes 1$ | $i \otimes i$ | $-1 \otimes j$ | $1 \otimes k$ | $1 \otimes 1$ | $-1 \otimes i$ |
| $i \otimes k$ | $i \otimes j$ | $-i \otimes i$ | $-i \otimes 1$ | $-1 \otimes k$ | $-1 \otimes j$ | $1 \otimes i$ | $1 \otimes 1$ |

*Proof.* Table is written according to equation

$$(\overline{e}_{C\cdot i} \otimes \overline{e}_{12\cdot k})(\overline{e}_{C\cdot j} \otimes \overline{e}_{12\cdot m}) = (\overline{e}_{C\cdot i}\,\overline{e}_{C\cdot j}) \otimes (\overline{e}_{12\cdot k}\,\overline{e}_{12\cdot m})$$

and definition of basis (6.2). □

**Theorem 6.2.** *Structural constants of the algebra $E(C, -1, -1)$ over field $R$ have form*

$$C_{13\cdot\cdot00\cdot00}^{\cdot\cdot00} = 1 \quad C_{13\cdot\cdot00\cdot01}^{\cdot\cdot01} = 1 \quad C_{13\cdot\cdot00\cdot02}^{\cdot\cdot02} = 1 \quad C_{13\cdot\cdot00\cdot03}^{\cdot\cdot03} = 1$$

$$C_{13\cdot\cdot01\cdot00}^{\cdot\cdot01} = 1 \quad C_{13\cdot\cdot01\cdot01}^{\cdot\cdot00} = -1 \quad C_{13\cdot\cdot01\cdot02}^{\cdot\cdot03} = 1 \quad C_{13\cdot\cdot01\cdot03}^{\cdot\cdot02} = -1$$

$$C_{13\cdot\cdot02\cdot00}^{\cdot\cdot02} = 1 \quad C_{13\cdot\cdot02\cdot01}^{\cdot\cdot03} = -1 \quad C_{13\cdot\cdot02\cdot02}^{\cdot\cdot00} = -1 \quad C_{13\cdot\cdot02\cdot03}^{\cdot\cdot01} = 1$$

$$C_{13\cdot\cdot03\cdot00}^{\cdot\cdot03} = 1 \quad C_{13\cdot\cdot03\cdot01}^{\cdot\cdot02} = 1 \quad C_{13\cdot\cdot03\cdot02}^{\cdot\cdot01} = -1 \quad C_{13\cdot\cdot03\cdot03}^{\cdot\cdot00} = -1$$

$$C_{13\cdot\cdot00\cdot10}^{\cdot\cdot10} = 1 \quad C_{13\cdot\cdot00\cdot11}^{\cdot\cdot11} = 1 \quad C_{13\cdot\cdot00\cdot12}^{\cdot\cdot12} = 1 \quad C_{13\cdot\cdot00\cdot13}^{\cdot\cdot13} = 1$$

$$C_{13\cdot\cdot01\cdot10}^{\cdot\cdot11} = 1 \quad C_{13\cdot\cdot01\cdot11}^{\cdot\cdot10} = -1 \quad C_{13\cdot\cdot01\cdot12}^{\cdot\cdot13} = 1 \quad C_{13\cdot\cdot01\cdot13}^{\cdot\cdot12} = -1$$

$$C_{13\cdot\cdot02\cdot10}^{\cdot\cdot12} = 1 \quad C_{13\cdot\cdot02\cdot11}^{\cdot\cdot13} = -1 \quad C_{13\cdot\cdot02\cdot12}^{\cdot\cdot10} = -1 \quad C_{13\cdot\cdot02\cdot13}^{\cdot\cdot11} = 1$$

$$C_{13\cdot\cdot03\cdot10}^{\cdot\cdot13} = 1 \quad C_{13\cdot\cdot03\cdot11}^{\cdot\cdot12} = 1 \quad C_{13\cdot\cdot03\cdot12}^{\cdot\cdot11} = -1 \quad C_{13\cdot\cdot03\cdot13}^{\cdot\cdot10} = -1$$

$$C_{13\cdot\cdot10\cdot00}^{\cdot\cdot10} = 1 \quad C_{13\cdot\cdot10\cdot01}^{\cdot\cdot11} = 1 \quad C_{13\cdot\cdot10\cdot02}^{\cdot\cdot12} = 1 \quad C_{13\cdot\cdot10\cdot03}^{\cdot\cdot13} = 1$$

$$C_{13\cdot\cdot11\cdot00}^{\cdot\cdot11} = 1 \quad C_{13\cdot\cdot11\cdot01}^{\cdot\cdot10} = -1 \quad C_{13\cdot\cdot11\cdot02}^{\cdot\cdot13} = 1 \quad C_{13\cdot\cdot11\cdot03}^{\cdot\cdot12} = -1$$

$$C_{13\cdot\cdot12\cdot00}^{\cdot\cdot12} = 1 \quad C_{13\cdot\cdot12\cdot01}^{\cdot\cdot13} = -1 \quad C_{13\cdot\cdot12\cdot02}^{\cdot\cdot10} = -1 \quad C_{13\cdot\cdot12\cdot03}^{\cdot\cdot11} = 1$$

$$C_{13\cdot\cdot13\cdot00}^{\cdot\cdot13} = 1 \quad C_{13\cdot\cdot13\cdot01}^{\cdot\cdot12} = 1 \quad C_{13\cdot\cdot13\cdot02}^{\cdot\cdot11} = -1 \quad C_{13\cdot\cdot13\cdot03}^{\cdot\cdot10} = -1$$

$$C_{13\cdot\cdot10\cdot10}^{\cdot\cdot00} = -1 \quad C_{13\cdot\cdot10\cdot11}^{\cdot\cdot01} = -1 \quad C_{13\cdot\cdot10\cdot12}^{\cdot\cdot02} = -1 \quad C_{13\cdot\cdot10\cdot13}^{\cdot\cdot03} = -1$$

$$C_{13\cdot\cdot11\cdot10}^{\cdot\cdot01} = -1 \quad C_{13\cdot\cdot11\cdot11}^{\cdot\cdot00} = 1 \quad C_{13\cdot\cdot11\cdot12}^{\cdot\cdot03} = -1 \quad C_{13\cdot\cdot11\cdot13}^{\cdot\cdot02} = 1$$

$$C_{13\cdot\cdot12\cdot10}^{\cdot\cdot02} = -1 \quad C_{13\cdot\cdot12\cdot11}^{\cdot\cdot03} = 1 \quad C_{13\cdot\cdot12\cdot12}^{\cdot\cdot00} = 1 \quad C_{13\cdot\cdot12\cdot13}^{\cdot\cdot01} = -1$$

$$C_{13\cdot\cdot13\cdot10}^{\cdot\cdot03} = -1 \quad C_{13\cdot\cdot13\cdot11}^{\cdot\cdot02} = -1 \quad C_{13\cdot\cdot13\cdot12}^{\cdot\cdot01} = 1 \quad C_{13\cdot\cdot13\cdot13}^{\cdot\cdot00} = 1$$

*Proof.* We consider the statement of theorem either as corollary of the theorem 6.1, or as corollary of theorem 5.1. □

**Theorem 6.3** (The Cauchy-Riemann equations). *Matrix of linear function*

$$y^{ik} = x^{jm} f_{jm}^{ik}$$

*of algebra $E(C, -1, -1)$ satisfies relationship*

(6.3)
$$f_{\cdot 0i}^{\cdot 0j} = f_{\cdot 1i}^{\cdot 1j}$$
$$f_{\cdot 0i}^{\cdot 1j} = -f_{\cdot 1i}^{\cdot 0j}$$

*Proof.* From equations [2]-(3.1.17), (5.8), (3.3) it follows

$$f_{\cdot 0i}^{\cdot 0j} = f^{\cdot ka\cdot rc} C_{13\cdot\cdot ka\cdot 0i}^{\cdot\cdot pb} C_{13\cdot\cdot pb\cdot rc}^{\cdot\cdot 0j}$$

$$= f^{\cdot ka\cdot rc} C_{C\cdot k0}^{\cdot p} C_{12\cdot ai}^{\cdot b} C_{C\cdot pr}^{\cdot 0} C_{12\cdot bc}^{\cdot j}$$

$$= f^{\cdot 0a\cdot 0c} C_{C\cdot 00}^{\cdot 0} C_{12\cdot ai}^{\cdot b} C_{C\cdot 00}^{\cdot 0} C_{12\cdot bc}^{\cdot j}$$

(6.4)
$$+ f^{\cdot 1a\cdot 1c} C_{C\cdot 10}^{\cdot 1} C_{12\cdot ai}^{\cdot b} C_{C\cdot 11}^{\cdot 0} C_{12\cdot bc}^{\cdot j}$$

$$= f^{\cdot 0a\cdot 0c} C_{12\cdot ai}^{\cdot b} C_{12\cdot bc}^{\cdot j} - f^{\cdot 1a\cdot 1c} C_{12\cdot ai}^{\cdot b} C_{12\cdot bc}^{\cdot j}$$

$$= (f^{\cdot 0a\cdot 0c} - f^{\cdot 1a\cdot 1c}) C_{12\cdot ai}^{\cdot b} C_{12\cdot bc}^{\cdot j}$$

$$f_{\cdot 0i}^{\cdot 1j} = f^{\cdot ka \cdot rc} \, C_{13 \cdot \cdot ka \cdot 0i}^{\cdot \cdot pb} \, C_{13 \cdot \cdot pb \cdot rc}^{\cdot 1j}$$
$$= f^{\cdot ka \cdot rc} \, C_{C \cdot k0}^{\cdot p} \, C_{12 \cdot ai}^{b} \, C_{C \cdot pr}^{\cdot 1} \, C_{12 \cdot bc}^{j}$$
$$= f^{\cdot 0a \cdot 1c} \, C_{C \cdot 00}^{\cdot 0} \, C_{12 \cdot ai}^{b} \, C_{C \cdot 01}^{\cdot 1} \, C_{12 \cdot bc}^{j}$$
(6.5)
$$+ f^{\cdot 1a \cdot 0c} \, C_{C \cdot 10}^{\cdot 1} \, C_{12 \cdot ai}^{b} \, C_{C \cdot 10}^{\cdot 1} \, C_{12 \cdot bc}^{j}$$
$$= f^{\cdot 0a \cdot 1c} \, C_{12 \cdot ai}^{b} \, C_{12 \cdot bc}^{j} + f^{\cdot 1a \cdot 0c} \, C_{12 \cdot ai}^{b} \, C_{12 \cdot bc}^{j}$$
$$= (f^{\cdot 0a \cdot 1c} + f^{\cdot 1a \cdot 0c}) \, C_{12 \cdot ai}^{b} \, C_{12 \cdot bc}^{j}$$

$$f_{\cdot 1i}^{\cdot 0j} = f^{\cdot ka \cdot rc} \, C_{13 \cdot \cdot ka \cdot 1i}^{\cdot \cdot pb} \, C_{13 \cdot \cdot pb \cdot rc}^{\cdot 0j}$$
$$= f^{\cdot ka \cdot rc} \, C_{C \cdot k1}^{\cdot p} \, C_{12 \cdot ai}^{b} \, C_{C \cdot pr}^{\cdot 0} \, C_{12 \cdot bc}^{j}$$
$$= f^{\cdot 0a \cdot 1c} \, C_{C \cdot 01}^{\cdot 1} \, C_{12 \cdot ai}^{b} \, C_{C \cdot 11}^{\cdot 0} \, C_{12 \cdot bc}^{j}$$
(6.6)
$$+ f^{\cdot 1a \cdot 0c} \, C_{C \cdot 11}^{\cdot 0} \, C_{12 \cdot ai}^{b} \, C_{C \cdot 00}^{\cdot 0} \, C_{12 \cdot bc}^{j}$$
$$= -f^{\cdot 0a \cdot 1c} \, C_{12 \cdot ai}^{b} \, C_{12 \cdot bc}^{j} - f^{\cdot 1a \cdot 0c} \, C_{12 \cdot ai}^{b} \, C_{12 \cdot bc}^{j}$$
$$= -(f^{\cdot 0a \cdot 1c} + f^{\cdot 1a \cdot 0c}) \, C_{12 \cdot ai}^{b} \, C_{12 \cdot bc}^{j}$$

$$f_{\cdot 1i}^{\cdot 1j} = f^{\cdot ka \cdot rc} \, C_{13 \cdot \cdot ka \cdot 1i}^{\cdot \cdot pb} \, C_{13 \cdot \cdot pb \cdot rc}^{\cdot 1j}$$
$$= f^{\cdot ka \cdot rc} \, C_{C \cdot k1}^{\cdot p} \, C_{12 \cdot ai}^{b} \, C_{C \cdot pr}^{\cdot 1} \, C_{12 \cdot bc}^{j}$$
$$= f^{\cdot 0a \cdot 0c} \, C_{C \cdot 01}^{\cdot 1} \, C_{12 \cdot ai}^{b} \, C_{C \cdot 10}^{\cdot 1} \, C_{12 \cdot bc}^{j}$$
(6.7)
$$+ f^{\cdot 1a \cdot 1c} \, C_{C \cdot 11}^{\cdot 0} \, C_{12 \cdot ai}^{b} \, C_{C \cdot 01}^{\cdot 1} \, C_{12 \cdot bc}^{j}$$
$$= f^{\cdot 0a \cdot 0c} \, C_{12 \cdot ai}^{b} \, C_{12 \cdot bc}^{j} - f^{\cdot 1a \cdot 1c} \, C_{12 \cdot ai}^{b} \, C_{12 \cdot bc}^{j}$$
$$= (f^{\cdot 0a \cdot 0c} - f^{\cdot 1a \cdot 1c}) \, C_{12 \cdot ai}^{b} \, C_{12 \cdot bc}^{j}$$

Equation (6.3) follows from comparison of equations (6.4) and (6.7), (6.5) and (6.6). □

**Theorem 6.4** (The Cauchy-Riemann equations). *Since matrix*

$$\left(\frac{\partial y^{\cdot ij}}{\partial x^{\cdot kl}}\right)$$

*is Jacobian matrix of map in algebra* $E(C, -1, -1)$, *then*

$$\frac{\partial y^{\cdot 1j}}{\partial x^{\cdot 0i}} = -\frac{\partial y^{\cdot 0j}}{\partial x^{\cdot 1i}}$$
$$\frac{\partial y^{\cdot 0j}}{\partial x^{\cdot 0i}} = \frac{\partial y^{\cdot 1j}}{\partial x^{\cdot 1i}}$$

*Proof.* The statement of theorem is corollary of theorem 6.3. □

## 7. Algebra $C \otimes (C \otimes H)$

Algebra $C \otimes (C \otimes H)$ is not quaternion algebra. However I consider this algebra here because this algebra is similar on algebra $E(C, -1, -1) = C \otimes H$.

Algebra $C \otimes (C \otimes H)$ is interesting from other point of view also. When we consider this algebra it becomes evident that tensor product of noncommutative rings is nonassociative. Because $C$ is field then $C \otimes C$ is isomorphic to $C$, and therefore, $(C \otimes C) \otimes H$ is isomorphic $C \otimes H$. However, because $H$ is not algebra over complex field, then algebra $C \otimes (C \otimes H)$ is different from algebra $(C \otimes C) \otimes H = C \otimes H$.

Algebra $C \otimes (C \otimes H)$ has so large dimension, that it becomes unsuitable to consider product table and structural constants of this algebra. However it is easy to see that structure of this algebra is similar to structure of algebra considered in section 5.

We will represent basis of algebra $C \otimes (C \otimes H)$ as

$$\overline{e}_{kji} = \overline{e}_{C \cdot k} \otimes (\overline{e}_{C \cdot j} \otimes \overline{e}_{H \cdot i}) \tag{7.1}$$

where $\overline{\overline{e}}_C$ is basis of algebra $C$ and $\overline{\overline{e}}_H$ is basis of algebra $H$. Correspondingly, the product in algebra $C \otimes (C \otimes H)$ is defined componentwise

$$a_1 \otimes (a_2 \otimes a_3)\, b_1 \otimes (b_2 \otimes b_2) = (a_1 b_1) \otimes ((a_2 b_2) \otimes (a_3 b_3)) \tag{7.2}$$

**Theorem 7.1.** *Structural constants of the algebra* $C \otimes (C \otimes H)$ *have form*

$$B_{\cdot pji \cdot qmk}^{\cdot rdb} = C_{C \cdot pq}^{\ r}\ C_{C \cdot jm}^{\ d}\ C_{H \cdot ik}^{\ b}$$

*where* $C_{C \cdot pq}^{\ r}$ *are structural constants of complex field,* $C_{H \cdot ik}^{\ b}$ *are structural constants of quaternion algebra.*

*Proof.* To prove the theorem it is enough to compare following equations

$$\begin{aligned}
\overline{e}_{pji}\, \overline{e}_{qmk} &= B_{\cdot pji \cdot qmk}^{\cdot rdb}\, \overline{e}_{rdb} \\
\overline{e}_{pji}\, \overline{e}_{qmk} &= \overline{e}_{C \cdot p} \otimes (\overline{e}_{C \cdot j} \otimes \overline{e}_{H \cdot i}) \overline{e}_{C \cdot q} \otimes (\overline{e}_{C \cdot m} \otimes \overline{e}_{H \cdot k})) \\
&= (\overline{e}_{C \cdot p}\, \overline{e}_{C \cdot q}) \otimes ((\overline{e}_{C \cdot j}\, \overline{e}_{C \cdot m}) \otimes (\overline{e}_{H \cdot i})\, \overline{e}_{H \cdot k})) \\
&= C_{C \cdot pq}^{\ r} C_{C \cdot jm}^{\ d} C_{H \cdot ik}^{\ b} \overline{e}_{C \cdot r} \otimes (\overline{e}_{C \cdot d} \otimes \overline{e}_{H \cdot b}) \\
&= C_{C \cdot pq}^{\ r}\ C_{C \cdot jm}^{\ d}\ C_{H \cdot ik}^{\ b}\ \overline{e}_{rdb}
\end{aligned}$$

$\square$

**Theorem 7.2.** *Structural constants of the algebra* $C \otimes (C \otimes H)$ *have form*

$$\begin{aligned}
B_{\cdot 00i \cdot 00k}^{\cdot 00b} &= \phantom{-}C_{H \cdot ik}^{\ b} & B_{\cdot 00i \cdot 01k}^{\cdot 01b} &= \phantom{-}C_{H \cdot ik}^{\ b} \\
B_{\cdot 01i \cdot 00k}^{\cdot 01b} &= \phantom{-}C_{H \cdot ik}^{\ b} & B_{\cdot 01i \cdot 01k}^{\cdot 01b} &= -C_{H \cdot ik}^{\ b} \\
B_{\cdot 00i \cdot 10k}^{\cdot 10b} &= \phantom{-}C_{H \cdot ik}^{\ b} & B_{\cdot 00i \cdot 11k}^{\cdot 11b} &= \phantom{-}C_{H \cdot ik}^{\ b} \\
B_{\cdot 01i \cdot 10k}^{\cdot 11b} &= \phantom{-}C_{H \cdot ik}^{\ b} & B_{\cdot 01i \cdot 11k}^{\cdot 11b} &= -C_{H \cdot ik}^{\ b} \\
B_{\cdot 10i \cdot 00k}^{\cdot 10b} &= \phantom{-}C_{H \cdot ik}^{\ b} & B_{\cdot 10i \cdot 01k}^{\cdot 11b} &= \phantom{-}C_{H \cdot ik}^{\ b} \\
B_{\cdot 11i \cdot 00k}^{\cdot 11b} &= \phantom{-}C_{H \cdot ik}^{\ b} & B_{\cdot 11i \cdot 01k}^{\cdot 11b} &= -C_{H \cdot ik}^{\ b} \\
B_{\cdot 10i \cdot 10k}^{\cdot 00b} &= -C_{H \cdot ik}^{\ b} & B_{\cdot 10i \cdot 11k}^{\cdot 01b} &= -C_{H \cdot ik}^{\ b} \\
B_{\cdot 11i \cdot 10k}^{\cdot 01b} &= -C_{H \cdot ik}^{\ b} & B_{\cdot 11i \cdot 11k}^{\cdot 01b} &= \phantom{-}C_{H \cdot ik}^{\ b}
\end{aligned} \tag{7.3}$$

*Proof.* The statement of theorem is corollary of theorems 7.1, 3.1. $\square$

**Theorem 7.3** (The Cauchy-Riemann equations)**.** *Matrix of linear function*

$$y^{ikp} = x^{jmr} f_{jmr}^{ikp}$$

of algebra $C \otimes (C \otimes H)$ satisfies relationship

$$
\begin{aligned}
f_{\cdot 00i}^{\cdot 00j} &= f_{\cdot 10i}^{\cdot 10j} = f_{\cdot 01i}^{\cdot 01j} = f_{\cdot 11i}^{\cdot 11j} \\
f_{\cdot 00i}^{\cdot 01j} &= f_{\cdot 10i}^{\cdot 11j} = -f_{\cdot 01i}^{\cdot 00j} = -f_{\cdot 11i}^{\cdot 10j} \\
f_{\cdot 00i}^{\cdot 10j} &= -f_{\cdot 10i}^{\cdot 00j} = f_{\cdot 01i}^{\cdot 11j} = -f_{\cdot 11i}^{\cdot 01j} \\
f_{\cdot 00i}^{\cdot 11j} &= -f_{\cdot 10i}^{\cdot 01j} = -f_{\cdot 01i}^{\cdot 10j} = f_{\cdot 11i}^{\cdot 00j}
\end{aligned}
\tag{7.4}
$$

*Proof.* Using equations [2]-(3.1.17), (7.3), we can repeat computing that we used to prove theorem 6.3. However, it is evident that this computing is different only by set of indexes. Thus, from theorem 6.3, it follows

$$
\begin{aligned}
f_{\cdot 0ji}^{\cdot 0mk} &= f_{\cdot 1ji}^{\cdot 1mk} \\
f_{\cdot 0ji}^{\cdot 1mk} &= -f_{\cdot 1ji}^{\cdot 0mk}
\end{aligned}
\tag{7.5}
$$

$$
\begin{aligned}
f_{\cdot p0i}^{\cdot q0k} &= f_{\cdot p1i}^{\cdot q1k} \\
f_{\cdot p0i}^{\cdot q1k} &= -f_{\cdot p1i}^{\cdot q0k}
\end{aligned}
\tag{7.6}
$$

Equations (7.4) follow from equations (7.5), (7.6). □

**Theorem 7.4** (The Cauchy-Riemann equations). *Since matrix*

$$\left(\frac{\partial y^{\cdot ijp}}{\partial x^{\cdot klq}}\right)$$

*is Jacobian matrix of map in algebra* $C \otimes (C \otimes H)$ *, then*

$$
\frac{\partial y^{\cdot 00i}}{\partial x^{\cdot 00j}} = \frac{\partial y^{\cdot 10i}}{\partial x^{\cdot 10j}} = \frac{\partial y^{\cdot 01i}}{\partial x^{\cdot 01j}} = \frac{\partial y^{\cdot 11i}}{\partial x^{\cdot 11j}}
$$

$$
\frac{\partial y^{\cdot 00i}}{\partial x^{\cdot 01j}} = \frac{\partial y^{\cdot 10i}}{\partial x^{\cdot 11j}} = -\frac{\partial y^{\cdot 01i}}{\partial x^{\cdot 00j}} = -\frac{\partial y^{\cdot 11i}}{\partial x^{\cdot 10j}}
$$

$$
\frac{\partial y^{\cdot 00i}}{\partial x^{\cdot 10j}} = -\frac{\partial y^{\cdot 10i}}{\partial x^{\cdot 00j}} = \frac{\partial y^{\cdot 01i}}{\partial x^{\cdot 11j}} = -\frac{\partial y^{\cdot 11i}}{\partial x^{\cdot 01j}}
$$

$$
\frac{\partial y^{\cdot 00i}}{\partial x^{\cdot 11j}} = -\frac{\partial y^{\cdot 10i}}{\partial x^{\cdot 01j}} = -\frac{\partial y^{\cdot 01i}}{\partial x^{\cdot 10j}} = \frac{\partial y^{\cdot 11i}}{\partial x^{\cdot 00j}}
$$

*Proof.* The statement of theorem is corollary of theorem 7.3. □

## 8. Quaternion Algebra $E(R, a, b)$

Assume $\overline{e}_0 = 1$, $\overline{e}_1 = i$, $\overline{e}_2 = j$, $\overline{e}_3 = k$. According to equation (4.1) structural constants of algebra $E(R, a, b)$ have form

$$
\begin{aligned}
C_{00}^0 &= 1 & C_{01}^1 &= 1 & C_{02}^2 &= 1 & C_{03}^3 &= 1 \\
C_{10}^1 &= 1 & C_{11}^0 &= a & C_{12}^3 &= 1 & C_{13}^2 &= a \\
C_{20}^2 &= 1 & C_{21}^3 &= -1 & C_{22}^0 &= b & C_{23}^1 &= -b \\
C_{30}^3 &= 1 & C_{31}^2 &= -a & C_{32}^1 &= b & C_{33}^0 &= -ab
\end{aligned}
$$

**Theorem 8.1.** *Standard components of linear function and coordinates of corresponding linear map over field $R$ satisfy relationship*

(8.1)
$$\begin{cases} f_0^0 = f^{00} + af^{11} + bf^{22} - abf^{33} \\ f_1^1 = f^{00} + af^{11} - bf^{22} + abf^{33} \\ f_2^2 = f^{00} - af^{11} + bf^{22} + abf^{33} \\ f_3^3 = f^{00} - af^{11} - bf^{22} - abf^{33} \end{cases}$$

(8.2)
$$\begin{cases} f_0^1 = f^{01} + f^{10} - bf^{23} + bf^{32} \\ f_1^0 = af^{01} + af^{10} + abf^{23} - abf^{32} \\ f_2^3 = -f^{01} + f^{10} + bf^{23} + bf^{32} \\ f_3^2 = -af^{01} + af^{10} - abf^{23} - abf^{32} \end{cases}$$

(8.3)
$$\begin{cases} f_0^2 = f^{02} + af^{13} + f^{20} - af^{31} \\ f_1^3 = f^{02} + af^{13} - f^{20} + af^{31} \\ f_2^0 = bf^{02} - abf^{13} + bf^{20} + abf^{31} \\ f_3^1 = bf^{02} - abf^{13} - bf^{20} - abf^{31} \end{cases}$$

(8.4)
$$\begin{cases} f_0^3 = f^{03} + f^{12} - f^{21} + f^{30} \\ f_1^2 = af^{03} + af^{12} + af^{21} - af^{30} \\ f_2^1 = -bf^{03} + bf^{12} + bf^{21} + bf^{30} \\ f_3^0 = -abf^{03} + abf^{12} - abf^{21} - abf^{30} \end{cases}$$

(8.5)
$$\begin{cases} f^{00} = \dfrac{1}{4} f_0^0 + \dfrac{1}{4} f_1^1 + \dfrac{1}{4} f_2^2 + \dfrac{1}{4} f_3^3 \\ f^{11} = \dfrac{1}{4a} f_0^0 + \dfrac{1}{4a} f_1^1 - \dfrac{1}{4a} f_2^2 - \dfrac{1}{4a} f_3^3 \\ f^{22} = \dfrac{1}{4b} f_0^0 - \dfrac{1}{4b} f_1^1 + \dfrac{1}{4b} f_2^2 - \dfrac{1}{4b} f_3^3 \\ f^{33} = -\dfrac{1}{4ab} f_0^0 + \dfrac{1}{4ab} f_1^1 + \dfrac{1}{4ab} f_2^2 - \dfrac{1}{4ab} f_3^3 \end{cases}$$

(8.6)
$$\begin{cases} f^{10} = \dfrac{1}{4a} f_1^0 + \dfrac{1}{4} f_0^1 + \dfrac{1}{4a} f_3^2 + \dfrac{1}{4} f_2^3 \\ f^{01} = \dfrac{1}{4a} f_1^0 + \dfrac{1}{4} f_0^1 - \dfrac{1}{4a} f_3^2 - \dfrac{1}{4} f_2^3 \\ f^{32} = -\dfrac{1}{4ab} f_1^0 + \dfrac{1}{4b} f_0^1 - \dfrac{1}{4ab} f_3^2 + \dfrac{1}{4b} f_2^3 \\ f^{23} = \dfrac{1}{4ab} f_1^0 - \dfrac{1}{4b} f_0^1 - \dfrac{1}{4ab} f_3^2 + \dfrac{1}{4b} f_2^3 \end{cases}$$

$$(8.7)\begin{cases} f^{20} = \phantom{-}\dfrac{1}{4b}\,f_2^0 - \dfrac{1}{4b}\,f_3^1 + \dfrac{1}{4}\,f_0^2 - \dfrac{1}{4}\,f_1^3 \\[4pt] f^{31} = \phantom{-}\dfrac{1}{4ab}f_2^0 - \dfrac{1}{4ab}f_3^1 - \dfrac{1}{4a}f_0^2 + \dfrac{1}{4a}f_1^3 \\[4pt] f^{02} = \phantom{-}\dfrac{1}{4b}\,f_2^0 + \dfrac{1}{4b}\,f_3^1 + \dfrac{1}{4}\,f_0^2 + \dfrac{1}{4}\,f_1^3 \\[4pt] f^{13} = -\dfrac{1}{4ab}f_2^0 - \dfrac{1}{4ab}f_3^1 + \dfrac{1}{4a}f_0^2 + \dfrac{1}{4a}f_1^3 \end{cases}$$

$$(8.8)\begin{cases} f^{30} = -\dfrac{1}{4ab}f_3^0 + \dfrac{1}{4b}f_2^1 - \dfrac{1}{4a}f_1^2 + \dfrac{1}{4}f_0^3 \\[4pt] f^{21} = -\dfrac{1}{4ab}f_3^0 + \dfrac{1}{4b}f_2^1 + \dfrac{1}{4a}f_1^2 - \dfrac{1}{4}f_0^3 \\[4pt] f^{12} = \phantom{-}\dfrac{1}{4ab}f_3^0 + \dfrac{1}{4b}f_2^1 + \dfrac{1}{4a}f_1^2 + \dfrac{1}{4}f_0^3 \\[4pt] f^{03} = -\dfrac{1}{4ab}f_3^0 - \dfrac{1}{4b}f_2^1 + \dfrac{1}{4a}f_1^2 + \dfrac{1}{4}f_0^3 \end{cases}$$

*Proof.* Using equation [2]-3.1.17 we get relationships

$$f_0^0 = f^{kr}\,C_{k0}^p\,C_{pr}^0$$
$$= f^{00}\,C_{00}^0\,C_{00}^0 + f^{11}\,C_{10}^1\,C_{11}^0 + f^{22}\,C_{20}^2\,C_{22}^0 + f^{33}\,C_{30}^3\,C_{33}^0$$
$$= f^{00} + af^{11} + bf^{22} - abf^{33}$$

$$f_0^1 = f^{kr}\,C_{k0}^p\,C_{pr}^1$$
$$= f^{01}\,C_{00}^0\,C_{01}^1 + f^{10}\,C_{10}^1\,C_{10}^1 + f^{23}\,C_{20}^2\,C_{23}^1 + f^{32}\,C_{30}^3\,C_{32}^1$$
$$= f^{01} + f^{10} - bf^{23} + bf^{32}$$

$$f_0^2 = f^{kr}\,C_{k0}^p\,C_{pr}^2$$
$$= f^{02}\,C_{00}^0\,C_{02}^2 + f^{13}\,C_{10}^1\,C_{13}^2 + f^{20}\,C_{20}^2\,C_{20}^2 + f^{31}\,C_{30}^3\,C_{31}^2$$
$$= f^{02} + af^{13} + f^{20} - af^{31}$$

$$f_0^3 = f^{kr}\,C_{k0}^p\,C_{pr}^3$$
$$= f^{03}\,C_{00}^0\,C_{03}^3 + f^{12}\,C_{10}^1\,C_{12}^3 + f^{21}\,C_{20}^2\,C_{21}^3 + f^{30}\,C_{30}^3\,C_{30}^3$$
$$= f^{03} + f^{12} - f^{21} + f^{30}$$

$$f_1^0 = f^{kr}\,C_{k1}^p\,C_{pr}^0$$
$$= f^{01}\,C_{01}^1\,C_{11}^0 + f^{10}\,C_{11}^0\,C_{00}^0 + f^{23}\,C_{21}^3\,C_{33}^0 + f^{32}\,C_{31}^2\,C_{22}^0$$
$$= af^{01} + af^{10} + abf^{23} - abf^{32}$$

$$f_1^1 = f^{kr}\,C_{k1}^p\,C_{pr}^1$$
$$= f^{00}\,C_{01}^1\,C_{10}^1 + f^{11}\,C_{11}^0\,C_{01}^1 + f^{22}\,C_{21}^3\,C_{32}^1 + f^{33}\,C_{31}^2\,C_{23}^1$$
$$= f^{00} + af^{11} - bf^{22} + abf^{33}$$

$$\begin{aligned}
f_1^2 &= f^{kr}\ C_{k1}^p\ C_{pr}^2 \\
&= f^{03}\ C_{01}^1\ C_{13}^2 + f^{12}\ C_{11}^0\ C_{02}^2 + f^{21}\ C_{21}^3\ C_{31}^2 + f^{30}\ C_{31}^2\ C_{20}^2 \\
&= af^{03} + af^{12} + af^{21} - af^{30}
\end{aligned}$$

$$\begin{aligned}
f_1^3 &= f^{kr}\ C_{k1}^p\ C_{pr}^3 \\
&= f^{02}\ C_{01}^1\ C_{12}^3 + f^{13}\ C_{11}^0\ C_{03}^3 + f^{20}\ C_{21}^3\ C_{30}^3 + f^{31}\ C_{31}^2\ C_{21}^3 \\
&= f^{02} + af^{13} - f^{20} + af^{31}
\end{aligned}$$

$$\begin{aligned}
f_2^0 &= f^{kr}\ C_{k2}^p\ C_{pr}^0 \\
&= f^{02}\ C_{02}^2\ C_{22}^0 + f^{13}\ C_{12}^3\ C_{33}^0 + f^{20}\ C_{22}^0\ C_{00}^0 + f^{31}\ C_{32}^1\ C_{11}^0 \\
&= bf^{02} - abf^{13} + bf^{20} + abf^{31}
\end{aligned}$$

$$\begin{aligned}
f_2^1 &= f^{kr}\ C_{k2}^p\ C_{pr}^1 \\
&= f^{03}\ C_{02}^2\ C_{23}^1 + f^{12}\ C_{12}^3\ C_{32}^1 + f^{21}\ C_{22}^0\ C_{01}^1 + f^{30}\ C_{32}^1\ C_{10}^1 \\
&= -bf^{03} + bf^{12} + bf^{21} + bf^{30}
\end{aligned}$$

$$\begin{aligned}
f_2^2 &= f^{kr}\ C_{k2}^p\ C_{pr}^2 \\
&= f^{00}\ C_{02}^2\ C_{20}^2 + f^{11}\ C_{12}^3\ C_{31}^2 + f^{22}\ C_{22}^0\ C_{02}^2 + f^{33}\ C_{32}^1\ C_{13}^2 \\
&= f^{00} - af^{11} + bf^{22} + abf^{33}
\end{aligned}$$

$$\begin{aligned}
f_2^3 &= f^{kr}\ C_{k2}^p\ C_{pr}^3 \\
&= f^{01}\ C_{02}^2\ C_{21}^3 + f^{10}\ C_{12}^3\ C_{30}^3 + f^{23}\ C_{22}^0\ C_{03}^3 + f^{32}\ C_{32}^1\ C_{12}^3 \\
&= -f^{01} + f^{10} + bf^{23} + bf^{32}
\end{aligned}$$

$$\begin{aligned}
f_3^0 &= f^{kr}\ C_{k3}^p\ C_{pr}^0 \\
&= f^{03}\ C_{03}^3\ C_{33}^0 + f^{12}\ C_{13}^2\ C_{22}^0 + f^{21}\ C_{23}^1\ C_{11}^0 + f^{30}\ C_{33}^0\ C_{00}^0 \\
&= -abf^{03} + abf^{12} - abf^{21} - abf^{30}
\end{aligned}$$

$$\begin{aligned}
f_3^1 &= f^{kr}\ C_{k3}^p\ C_{pr}^1 \\
&= f^{02}\ C_{03}^3\ C_{32}^1 + f^{13}\ C_{13}^2\ C_{23}^1 + f^{20}\ C_{23}^1\ C_{10}^1 + f^{31}\ C_{33}^0\ C_{01}^1 \\
&= bf^{02} - abf^{13} - bf^{20} - abf^{31}
\end{aligned}$$

$$\begin{aligned}
f_3^2 &= f^{kr}\ C_{k3}^p\ C_{pr}^2 \\
&= f^{01}\ C_{03}^3\ C_{31}^2 + f^{10}\ C_{13}^2\ C_{20}^2 + f^{23}\ C_{23}^1\ C_{13}^2 + f^{32}\ C_{33}^0\ C_{02}^2 \\
&= -af^{01} + af^{10} - abf^{23} - abf^{32}
\end{aligned}$$

$$\begin{aligned}
f_3^3 &= f^{kr}\ C_{k3}^p\ C_{pr}^3 \\
&= f^{00}\ C_{03}^3\ C_{30}^3 + f^{11}\ C_{13}^2\ C_{21}^3 + f^{22}\ C_{23}^1\ C_{12}^3 + f^{33}\ C_{33}^0\ C_{03}^3 \\
&= f^{00} - af^{11} - bf^{22} - abf^{33}
\end{aligned}$$

We group these relationships into systems of linear equations (8.1), (8.2), (8.3), (8.4).

(8.5) is solution of system of linear equations (8.1).
(8.6) is solution of system of linear equations (8.2).
(8.7) is solution of system of linear equations (8.3).
(8.8) is solution of system of linear equations (8.4). □

**Theorem 8.2.** *For any values of parameters $a \neq 0$, $b \neq 0$, there exists one to one map between coordinates of linear function of algebra $E(R, a, b)$ and its standard components.*

*Proof.* The statement of theorem is corollary of theorem 8.1. □

## 9. Regular Function

Although there is no analogue of the Cauchy-Riemann equations in quaternion algebra, in different papers mathematicians explore different sets of functions that have properties similar to properties of functions of complex variable. In [4, 5], there is definition of regular function that satisfies to equation

$$(9.1) \qquad \frac{\partial f}{\partial x^0} + i\frac{\partial f}{\partial x^1} + j\frac{\partial f}{\partial x^2} + k\frac{\partial f}{\partial x^3} = 0$$

**Theorem 9.1.** *Differential equation (9.1) is equivalent to system of differential equations*

$$(9.2) \qquad \begin{cases} \dfrac{\partial f^0}{\partial x^0} - \dfrac{\partial f^1}{\partial x^1} - \dfrac{\partial f^2}{\partial x^2} - \dfrac{\partial f^3}{\partial x^3} = 0 \\ \dfrac{\partial f^0}{\partial x^1} + \dfrac{\partial f^1}{\partial x^0} - \dfrac{\partial f^2}{\partial x^3} + \dfrac{\partial f^3}{\partial x^2} = 0 \\ \dfrac{\partial f^0}{\partial x^2} + \dfrac{\partial f^1}{\partial x^3} + \dfrac{\partial f^2}{\partial x^0} - \dfrac{\partial f^3}{\partial x^1} = 0 \\ \dfrac{\partial f^0}{\partial x^3} + \dfrac{\partial f^2}{\partial x^1} - \dfrac{\partial f^1}{\partial x^2} + \dfrac{\partial f^3}{\partial x^0} = 0 \end{cases}$$

*Proof.* Let us substitute

$$(9.3) \qquad \frac{\partial f}{\partial x^i} = \frac{\partial f^0}{\partial x^i} + \frac{\partial f^1}{\partial x^i}i + \frac{\partial f^2}{\partial x^i}j + \frac{\partial f^3}{\partial x^i}k$$

into equation (9.1). We will get

$$(9.4) \qquad \begin{aligned} \frac{\partial f}{\partial x^0} + i\frac{\partial f}{\partial x^1} + j\frac{\partial f}{\partial x^2} + k\frac{\partial f}{\partial x^3} =& \frac{\partial f^0}{\partial x^0} - \frac{\partial f^1}{\partial x^1} - \frac{\partial f^2}{\partial x^2} - \frac{\partial f^3}{\partial x^3} \\ &+ i(\frac{\partial f^0}{\partial x^1} + \frac{\partial f^1}{\partial x^0} - \frac{\partial f^2}{\partial x^3} + \frac{\partial f^3}{\partial x^2}) \\ &+ j(\frac{\partial f^0}{\partial x^2} + \frac{\partial f^1}{\partial x^3} + \frac{\partial f^2}{\partial x^0} - \frac{\partial f^3}{\partial x^1}) \\ &+ k(\frac{\partial f^0}{\partial x^3} + \frac{\partial f^2}{\partial x^1} - \frac{\partial f^1}{\partial x^2} + \frac{\partial f^3}{\partial x^0}) = 0 \end{aligned}$$

The statement of theorem follows from equation (9.4). □

In the paper [1], corollary 3.1.2, p. 1000, Deavours proves that the only regular quaternion functions with bounded norm is a constant.

**Theorem 9.2.** *Components of the Gâteaux derivative of regular quaternion function satisfy to equations*

(9.5)
$$\frac{\partial^{00}y}{\partial x} + \frac{\partial^{11}y}{\partial x} + \frac{\partial^{22}y}{\partial x} + \frac{\partial^{33}y}{\partial x} = 0$$
$$-\frac{\partial^{01}y}{\partial x} + \frac{\partial^{10}y}{\partial x} + \frac{\partial^{23}y}{\partial x} - \frac{\partial^{32}y}{\partial x} = 0$$
$$-\frac{\partial^{02}y}{\partial x} - \frac{\partial^{13}y}{\partial x} + \frac{\partial^{20}y}{\partial x} + \frac{\partial^{31}y}{\partial x} = 0$$
$$-\frac{\partial^{03}y}{\partial x} - \frac{\partial^{12}y}{\partial x} + \frac{\partial^{21}y}{\partial x} + \frac{\partial^{30}y}{\partial x} = 0$$

*Proof.* Let us substitute equations [2]-(7.4.1), [2]-(7.4.2), [2]-(7.4.3), [2]-(7.4.4) into equation (9.2). We will get

$$\frac{\partial f^0}{\partial x^0} - \frac{\partial f^1}{\partial x^1} - \frac{\partial f^2}{\partial x^2} - \frac{\partial f^3}{\partial x^3}$$
$$= \frac{\partial^{00}y}{\partial x} - \frac{\partial^{11}y}{\partial x} - \frac{\partial^{22}y}{\partial x} - \frac{\partial^{33}y}{\partial x} - \left(\frac{\partial^{00}y}{\partial x} - \frac{\partial^{11}y}{\partial x} + \frac{\partial^{22}y}{\partial x} + \frac{\partial^{33}y}{\partial x}\right)$$
$$- \left(\frac{\partial^{00}y}{\partial x} + \frac{\partial^{11}y}{\partial x} - \frac{\partial^{22}y}{\partial x} + \frac{\partial^{33}y}{\partial x}\right) - \left(\frac{\partial^{00}y}{\partial x} + \frac{\partial^{11}y}{\partial x} + \frac{\partial^{22}y}{\partial x} - \frac{\partial^{33}y}{\partial x}\right)$$
$$= -2\frac{\partial^{00}y}{\partial x} - 2\frac{\partial^{11}y}{\partial x} - 2\frac{\partial^{22}y}{\partial x} - 2\frac{\partial^{33}y}{\partial x}$$
$$= 0$$

$$\frac{\partial f^0}{\partial x^1} + \frac{\partial f^1}{\partial x^0} - \frac{\partial f^2}{\partial x^3} + \frac{\partial f^3}{\partial x^2}$$
$$= -\frac{\partial^{01}y}{\partial x} - \frac{\partial^{10}y}{\partial x} + \frac{\partial^{23}y}{\partial x} - \frac{\partial^{32}y}{\partial x} + \left(\frac{\partial^{01}y}{\partial x} + \frac{\partial^{10}y}{\partial x} + \frac{\partial^{23}y}{\partial x} - \frac{\partial^{32}y}{\partial x}\right)$$
$$- \left(\frac{\partial^{01}y}{\partial x} - \frac{\partial^{10}y}{\partial x} - \frac{\partial^{23}y}{\partial x} - \frac{\partial^{32}y}{\partial x}\right) + \left(-\frac{\partial^{01}y}{\partial x} + \frac{\partial^{10}y}{\partial x} - \frac{\partial^{23}y}{\partial x} - \frac{\partial^{32}y}{\partial x}\right)$$
$$= -2\frac{\partial^{01}y}{\partial x} + 2\frac{\partial^{10}y}{\partial x} + 2\frac{\partial^{23}y}{\partial x} - 2\frac{\partial^{32}y}{\partial x}$$
$$= 0$$

$$\frac{\partial f^0}{\partial x^2} + \frac{\partial f^1}{\partial x^3} + \frac{\partial f^2}{\partial x^0} - \frac{\partial f^3}{\partial x^1}$$
$$= -\frac{\partial^{02}y}{\partial x} - \frac{\partial^{13}y}{\partial x} - \frac{\partial^{20}y}{\partial x} + \frac{\partial^{31}y}{\partial x} + \left(-\frac{\partial^{02}y}{\partial x} - \frac{\partial^{13}y}{\partial x} + \frac{\partial^{20}y}{\partial x} - \frac{\partial^{31}y}{\partial x}\right)$$
$$+ \left(\frac{\partial^{02}y}{\partial x} - \frac{\partial^{13}y}{\partial x} + \frac{\partial^{20}y}{\partial x} + \frac{\partial^{31}y}{\partial x}\right) - \left(\frac{\partial^{02}y}{\partial x} - \frac{\partial^{13}y}{\partial x} - \frac{\partial^{20}y}{\partial x} - \frac{\partial^{31}y}{\partial x}\right)$$
$$= -2\frac{\partial^{02}y}{\partial x} - 2\frac{\partial^{13}y}{\partial x} + 2\frac{\partial^{20}y}{\partial x} + 2\frac{\partial^{31}y}{\partial x}$$
$$= 0$$

$$\frac{\partial f^0}{\partial x^3} + \frac{\partial f^2}{\partial x^1} - \frac{\partial f^1}{\partial x^2} + \frac{\partial f^3}{\partial x^0}$$
$$= -\frac{\partial^{03} y}{\partial x} + \frac{\partial^{12} y}{\partial x} - \frac{\partial^{21} y}{\partial x} - \frac{\partial^{30} y}{\partial x} - \left(-\frac{\partial^{03} y}{\partial x} - \frac{\partial^{12} y}{\partial x} - \frac{\partial^{21} y}{\partial x} + \frac{\partial^{30} y}{\partial x}\right)$$
$$- \left(\frac{\partial^{03} y}{\partial x} - \frac{\partial^{12} y}{\partial x} - \frac{\partial^{21} y}{\partial x} - \frac{\partial^{30} y}{\partial x}\right) - \left(\frac{\partial^{03} y}{\partial x} + \frac{\partial^{12} y}{\partial x} - \frac{\partial^{21} y}{\partial x} + \frac{\partial^{30} y}{\partial x}\right)$$
$$= -2\frac{\partial^{03} y}{\partial x} - 2\frac{\partial^{12} y}{\partial x} + 2\frac{\partial^{21} y}{\partial x} + 2\frac{\partial^{30} y}{\partial x}$$
$$=0$$

$\square$

**Theorem 9.3.** *The Gâteaux differential of regular function over quaternion algebra has form*

$$
\begin{aligned}
(9.6) \quad & -\left(\frac{\partial^{11} y}{\partial x} + \frac{\partial^{22} y}{\partial x} + \frac{\partial^{33} y}{\partial x}\right) dx + \frac{\partial^{11} y}{\partial x} i dx i + \frac{\partial^{22} y}{\partial x} j dx j + \frac{\partial^{33} y}{\partial x} k dx k \\
& + \left(\frac{\partial^{10} y}{\partial x} + \frac{\partial^{23} y}{\partial x} - \frac{\partial^{32} y}{\partial x}\right) dx i + \frac{\partial^{10} y}{\partial x} i dx + \frac{\partial^{23} y}{\partial x} j dx k - \frac{\partial^{32} y}{\partial x} k dx i \\
& + \left(-\frac{\partial^{13} y}{\partial x} + \frac{\partial^{20} y}{\partial x} + \frac{\partial^{31} y}{\partial x}\right) dx j - \frac{\partial^{13} y}{\partial x} i dx k + \frac{\partial^{20} y}{\partial x} j dx + \frac{\partial^{31} y}{\partial x} k dx i \\
& + \left(-\frac{\partial^{12} y}{\partial x} + \frac{\partial^{21} y}{\partial x} + \frac{\partial^{30} y}{\partial x}\right) dx k - \frac{\partial^{12} y}{\partial x} i dx j + \frac{\partial^{21} y}{\partial x} j dx i + \frac{\partial^{30} y}{\partial x} k dx
\end{aligned}
$$

*Proof.* The statement of theorem is corollary of theorem [2]-5.2.8. $\square$

Equation (9.1) is equivalent to equation

$$(9.7) \quad \begin{pmatrix} 1 & i & j & k \end{pmatrix} *^{*} \begin{pmatrix} \dfrac{\partial f^0}{\partial x^0} & \dfrac{\partial f^0}{\partial x^1} & \dfrac{\partial f^0}{\partial x^2} & \dfrac{\partial f^0}{\partial x^3} \\ \dfrac{\partial f^1}{\partial x^0} & \dfrac{\partial f^1}{\partial x^1} & \dfrac{\partial f^1}{\partial x^2} & \dfrac{\partial f^1}{\partial x^3} \\ \dfrac{\partial f^2}{\partial x^0} & \dfrac{\partial f^2}{\partial x^1} & \dfrac{\partial f^2}{\partial x^2} & \dfrac{\partial f^2}{\partial x^3} \\ \dfrac{\partial f^3}{\partial x^0} & \dfrac{\partial f^3}{\partial x^1} & \dfrac{\partial f^3}{\partial x^2} & \dfrac{\partial f^3}{\partial x^3} \end{pmatrix} *^{*} \begin{pmatrix} 1 \\ i \\ j \\ k \end{pmatrix} = 0$$

## 10. Instead of an Epilogue

The complex field and quaternion algebra have both common properties and differences. These differences make it harder to identify in the quaternion algebra patterns similar to those we have observed in the complex field. Therefore, it is very important to understand these differences.

One of the research directions is finding an analog to the Cauchy-Riemann equation in quaternion algebra. In this paper, I reviewed some studies in this area.

According to the theorem [2]-7.1.1, linear mapping has matrix

$$\begin{pmatrix} a_0 & -a_1 \\ a_1 & a_0 \end{pmatrix}$$

This mapping corresponds to multiplication by the number $a = a_0 + a_1 i$. The statement follows from equations

$$(a_0 + a_1 i)(x_0 + x_1 i) = a_0 x_0 - a_1 x_1 + (a_0 x_1 + a_1 x_0)i$$

$$\begin{pmatrix} a_0 & -a_1 \\ a_1 & a_0 \end{pmatrix} \begin{pmatrix} x_0 \\ x_1 \end{pmatrix} = \begin{pmatrix} a_0 x_0 - a_1 x_1 \\ a_1 x_0 + a_0 x_1 \end{pmatrix}$$

I decided to consider a similar class of functions of quaternions. The linear mapping of quaternion algebra

$$x \to ax$$

has matrix

(10.1) $$\begin{pmatrix} a^0 & -a^1 & -a^2 & -a^3 \\ a^1 & a^0 & -a^3 & a^2 \\ a^2 & a^3 & a^0 & -a^1 \\ a^3 & -a^2 & a^1 & a^0 \end{pmatrix}$$

It is interesting to consider the class of quaternion functions which has derivative of similar structure. However I think that the structure of the matrix (10.1) is too restrictive for derivative of functions of quaternions and I little relaxed the requirement. I assumed that the derivative satisfies to following equations

(10.2) $$\frac{\partial y^0}{\partial x^0} = \frac{\partial y^1}{\partial x^1} = \frac{\partial y^2}{\partial x^2} = \frac{\partial y^3}{\partial x^3}$$

(10.3) $$\frac{\partial y^i}{\partial x^j} = -\frac{\partial y^j}{\partial x^i} \quad i \neq j$$

It is easy to see that derivative of function like

$$y = ax \quad y = xa$$

satisfies to equations (10.2), (10.3).

Consider the function

$$y = x^2$$

Direct calculation gives

(10.4) $$\begin{aligned} y^0 &= (x^0)^2 - (x^1)^2 - (x^2)^2 - (x^3)^2 \\ y^1 &= 2x^0 x^1 \\ y^2 &= 2x^0 x^2 \\ y^3 &= 2x^0 x^3 \end{aligned}$$

The derivative of the mapping (10.4) has matrix

(10.5) $$\begin{pmatrix} 2x^0 & -2x^1 & -2x^2 & -2x^3 \\ 2x^1 & 2x^0 & 0 & 0 \\ 2x^2 & 0 & 2x^0 & 0 \\ 2x^3 & 0 & 0 & 2x^0 \end{pmatrix}$$

Therefore, the matrix (10.5) satisfies to equations (10.2), (10.3).

**Theorem 10.1.** *Derivative of conjugation does not satisfy equations* (10.2), (10.3).

*Proof.* According to the proof of the theorem [2]-7.3.5, derivative of conjugation does not satisfy equation (10.2). □

However, the set of functions whose derivative satisfies to equations (10.2), (10.3) is not sufficiently large.

**Theorem 10.2.** *The mapping*

$$y = x^3 \tag{10.6}$$

*does not satisfy to equation* (10.2).

*Proof.* We can represent the mapping (10.6) in the following form

$$y = x\, x^2$$

This mapping has coordinates

$$\begin{aligned}
y^0 &= x^0((x^0)^2 - (x^1)^2 - (x^2)^2 - (x^3)^2) \\
&\quad - x^1 2x^0 x^1 - x^2 2x^0 x^2 - x^3 2x^0 x^3 \\
&= (x^0)^3 - x^0(x^1)^2 - x^0(x^2)^2 - x^0(x^3)^2 \\
&\quad - 2x^0(x^1)^2 - 2x^0(x^2)^2 - 2x^0(x^3)^2 \\
&= (x^0)^3 - 3x^0(x^1)^2 - 3x^0(x^2)^2 - 3x^0(x^3)^2
\end{aligned} \tag{10.7}$$

$$\begin{aligned}
y^1 &= x^1((x^0)^2 - (x^1)^2 - (x^2)^2 - (x^3)^2) \\
&\quad + x^0 2x^0 x^1 - x^3 2x^0 x^2 + x^2 2x^0 x^3 \\
&= x^1(x^0)^2 - (x^1)^3 - x^1(x^2)^2 - x^1(x^3)^2 \\
&\quad + 2(x^0)^2 x^1 - 2x^0 x^2 x^3 + 2x^0 x^2 x^3 \\
&= 3x^1(x^0)^2 - (x^1)^3 - x^1(x^2)^2 - x^1(x^3)^2
\end{aligned} \tag{10.8}$$

From the equation (10.7), it follows that

$$\frac{\partial y^0}{\partial x^0} = 3(x^0)^2 - 3(x^1)^2 - 3(x^2)^2 - 3(x^3)^2 \tag{10.9}$$

From the equation (10.8), it follows that

$$\frac{\partial y^1}{\partial x^1} = 3(x^0)^2 - 3(x^1)^2 - (x^2)^2 - (x^3)^2 \tag{10.10}$$

The statement of the theorem follows from equations (10.9), (10.10). □

Without a doubt, this is only the beginning of the study, and many questions must be answered.

# Этюд о кватернионах

Александр Клейн

*Светлой памяти И. М. Гельфанда*

АННОТАЦИЯ. В статье рассматривается множество алгебр кватернионов над полем. Алгебра кватернионов $E(C, -1, -1)$ изоморфна тензорному произведению поля комплексных чисел $C$ и алгебры кватернионов $H = E(R, -1, -1)$. Рассмотренно множество функций алгебры кватернионов, которые могут иметь уравнения похожие на уравнение Коши-Римана для функций комплексной переменной.

## Содержание

## 1. Предисловие

Когда я начал исследовать математический анализ над телом, я обратил внимание на большое количество статей, посвящённое регулярным функциям кватернионов и аналогу уравнения Коши-Римана. Я хотел понять, появляется ли уравнение Коши-Римана в рамках той теории, которую я изучаю. Это и послужило толчком к рассмотрению тела как векторного пространства над центром.

Введение базиса упрощает некоторые построения и является своеобразным мостом между производной Гато и матрицей Якоби отображения. Это именно то место, где должно появиться уравнение Коши-Римана. Рассмотрение производной функции комплексных чисел показывает, что уравнение Коши-Римана имеет алгебраическую природу и связана с тем, что существует $R$-линейная функция поля комплексных чисел, не являющаяся $C$-линейной. Например сопряжение комплексного числа линейно над полем действительных чисел, но не является линейным отображением над полем комплексных чисел.

В алгебре кватернионов $H$ существует только $R$-линейное отображение. Следствием этого является возможность выразить любое $R$-линейное отображение через кватернион и отсутствие очевидного аналога уравнения Коши-Римана в алгебре кватернионов. Тем не менее, в многочисленных статьях и книгах, посвящённых анализу над алгеброй кватернионов, математики изучают множества отображений, свойства которых похожи на свойства функции комплексной переменной. Некоторые авторы не ограничивают себя алгеброй кватернионов и рассматривают более общие алгебры.

В статье [3] Гельфанд рассматривает алгебру кватернионов над произвольным полем и с произвольными параметрами. При этом $H = E(R, -1, -1)$. Опираясь на эту статью, я решил рассмотреть два случая, которые для меня важны.

Алгебра $E(R, a, b)$ была интересна для меня тем, что я хотел подобрать параметры $a$, $b$ так, что система линейных уравнений [2]-(3.1.17) становится вырожденной. Мне было важно понять, что произойдёт в этом случае. При изучении структуры линейного отображения над телом для меня было не совсем очевидно как повлияет на ответ вырожденность системы линейных уравнений [2]-(3.1.17). Однако решение этой задачи не дало ожидаемого ответа. Система линейных уравнений оказалась настолько простой, что невооружённым глазом видно, что эта система не может быть вырожденной.

Я писал, что уравнение Коши-Римана связано с тем, что поле комплексных чисел имеет подполе действительных чисел. Я высказал также предположение, что похожее утверждение возможно в алгебрах с достаточно сложным центром. Поэтому в алгебре кватернионов над полем комплексных чисел я ожидал увидеть аналог уравнения Коши-Римана. В ходе решения задачи выяснилось, что алгебра $E(C, -1, -1)$ изоморфна тензорному произведению $C \otimes H$. Поэтому линейные функции этой алгебры, а следовательно, и матрица Якоби произвольного отображения, удовлетворяют уравнению Коши-Римана для компоненты тензорного произведения, порождённой полем $C$. Естественным продолжением этой темы явилось изучение тензорного произведения $C \otimes (C \otimes H)$. В частности, первое, на что я обратил внимание, - это неассоциативность тензорного произведения.

Алгебра $C \otimes H$ интересна для меня и с другой точки зрения. Поле комплексных чисел является алгебраическим замыканием поля действительных чисел. Не возникает ли подобная связь между алгеброй $H$ и алгеброй $C \otimes H$?

В книге [2], на основе которой написана эта статья я рассматриваю производные Гато функции тела. Однако в этой статье я рассматриваю произвольные алгебры, которые не всегда являются телом. В течении того времени, что я работаю с производной Гато, я понял, что эта тема может быть обобщена на более широкий класс алгебр. Полное исследование появится позже, я написал эту статью с целью изучить условия появления уравнения Коши-Римана.

## 2. Соглашения

(1) Функция и отображение - синонимы. Однако существует традиция соответствие между кольцами или векторными пространствами называть отображением, а отображение поля действительных чисел или алгебры кватернионов называть функцией. Я тоже следую этой традиции, хотя встречается текст, в котором неясно, какому термину надо отдать предпочтение.

(2) Тело $D$ можно рассматривать как $D$-векторное пространство размерности 1. Соответственно этому, мы можем изучать не только гомоморфизм тела $D_1$ в тело $D_2$, но и линейное отображение тел. При этом подразумевается, что отображение мультипликативно над максимально возможным полем. В частности, линейное отображение тела $D$ мультипликативно над центром $Z(D)$. Это не противоречит определению линейного отображения поля, так как для поля $F$ справедливо $Z(F) = F$. Если поле $F$ отлично от максимально возможного, то я это явно указываю в тексте.

(3) Пусть $A$ - свободная конечно мерная алгебра. При разложении элемента алгебры $A$ относительно базиса $\overline{\overline{e}}$ мы пользуемся одной и той же корневой буквой для обозначения этого элемента и его координат. Однако в алгебре не принято использовать векторные обозначения. В выражении $a^2$ не ясно - это компонента разложения элемента $a$ относительно базиса или это операция возведения в степень. Для облегчения чтения текста мы будем индекс элемента алгебры выделять цветом. Например,
$$a = a^{i}\overline{e}_{i}$$

(4) Если свободная конечномерная алгебра имеет единицу, то мы будем отождествлять вектор базиса $\overline{e}_{0}$ с единицей алгебры.

(5) Без сомнения, у читателя могут быть вопросы, замечания, возражения. Я буду признателен любому отзыву.

## 3. Линейная функция комплексного поля

**Теорема 3.1** (Уравнения Коши-Римана). *Рассмотрим поле комплексных чисел $C$ как двумерную алгебру над полем действительных чисел. Положим*

$$\overline{e}_{C \cdot \mathbf{0}} = 1 \quad \overline{e}_{C \cdot \mathbf{1}} = i \tag{3.1}$$

*базис алгебры $C$. Тогда в этом базисе произведение имеет вид*

$$\overline{e}_{C \cdot \mathbf{1}}^{2} = -\overline{e}_{C \cdot \mathbf{0}} \tag{3.2}$$

*и структурные константы имеют вид*

$$\begin{aligned} C_{C \cdot \mathbf{00}}^{\mathbf{0}} &= 1 & C_{C \cdot \mathbf{01}}^{\mathbf{1}} &= 1 \\ C_{C \cdot \mathbf{10}}^{\mathbf{1}} &= 1 & C_{C \cdot \mathbf{11}}^{\mathbf{0}} &= -1 \end{aligned} \tag{3.3}$$

*Матрица линейной функции*
$$y^{\mathbf{i}} = x^{\mathbf{j}} f_{\mathbf{j}}^{\mathbf{i}}$$

поля комплексных чисел над полем действительных чисел удовлетворяет соотношению

$$f_0^0 = f_1^1 \tag{3.4}$$

$$f_0^1 = -f_1^0 \tag{3.5}$$

*Доказательство.* Равенства (3.2) и (3.3) следуют из равенства $i^2 = -1$. Пользуясь равенством [2]-(3.1.17) получаем соотношения

$$f_0^0 = f^{kr} C_{C \cdot k0}^{\ p} C_{C \cdot pr}^{\ 0} = f^{0r} C_{C \cdot 00}^{\ 0} C_{C \cdot 0r}^{\ 0} + f^{1r} C_{C \cdot 10}^{\ 1} C_{C \cdot 1r}^{\ 0} = f^{00} - f^{11} \tag{3.6}$$

$$f_0^1 = f^{kr} C_{C \cdot k0}^{\ p} C_{C \cdot pr}^{\ 1} = f^{0r} C_{C \cdot 00}^{\ 0} C_{C \cdot 0r}^{\ 1} + f^{1r} C_{C \cdot 10}^{\ 1} C_{C \cdot 1r}^{\ 1} = f^{01} + f^{10} \tag{3.7}$$

$$f_1^0 = f^{kr} C_{C \cdot k1}^{\ p} C_{C \cdot pr}^{\ 0} = f^{0r} C_{C \cdot 01}^{\ 1} C_{C \cdot 1r}^{\ 0} + f^{1r} C_{C \cdot 11}^{\ 0} C_{C \cdot 0r}^{\ 0} = -f^{01} - f^{10} \tag{3.8}$$

$$f_1^1 = f^{kr} C_{C \cdot k1}^{\ p} C_{C \cdot pr}^{\ 1} = f^{0r} C_{C \cdot 01}^{\ 1} C_{C \cdot 1r}^{\ 1} + f^{1r} C_{C \cdot 11}^{\ 0} C_{C \cdot 0r}^{\ 1} = f^{00} - f^{11} \tag{3.9}$$

Из равенств (3.6) и (3.9) следует (3.4). Из равенств (3.7) и (3.8) следует (3.5). □

**Теорема 3.2** (Уравнения Коши-Римана). *Если матрица*

$$\begin{pmatrix} \dfrac{\partial y^0}{\partial x^0} & \dfrac{\partial y^0}{\partial x^1} \\ \dfrac{\partial y^1}{\partial x^0} & \dfrac{\partial y^1}{\partial x^1} \end{pmatrix}$$

*является матрицей Якоби функции комплексного переменного*

$$x = x^0 + x^1 i \to y = y^0(x^0, x^1) + y^1(x^0, x^1)i$$

*над полем действительных чисел, то*

$$\begin{aligned} \frac{\partial y^1}{\partial x^0} &= -\frac{\partial y^0}{\partial x^1} \\ \frac{\partial y^0}{\partial x^0} &= \phantom{-}\frac{\partial y^1}{\partial x^1} \end{aligned} \tag{3.10}$$

*Доказательство.* Следствие теоремы 3.1. □

**Теорема 3.3.** *Производная функции комплексного переменного удовлетворяет равенству*

$$\frac{\partial y}{\partial x^0} + i\frac{\partial y}{\partial x^1} = 0 \tag{3.11}$$

*Доказательство.* Равенство

$$\frac{\partial y^0}{\partial x^0} + i\frac{\partial y^1}{\partial x^0} + i\frac{\partial y^0}{\partial x^1} - \frac{\partial y^1}{\partial x^1} = 0$$

следует из равенств (3.10). □

Равенство (3.11) эквивалентно равенству

$$\begin{pmatrix} 1 & i \end{pmatrix} \begin{pmatrix} \dfrac{\partial y^0}{\partial x^0} & \dfrac{\partial y^0}{\partial x^1} \\ \dfrac{\partial y^1}{\partial x^0} & \dfrac{\partial y^1}{\partial x^1} \end{pmatrix} \begin{pmatrix} 1 \\ i \end{pmatrix} = 0 \tag{3.12}$$

## 4. Алгебра кватернионов

В этой статье я рассматриваю множество кватернионных алгебр, определённых в [3].

**Определение 4.1.** Пусть $F$ - поле. Расширение $F(i,j,k)$ поля $F$ называется **алгеброй** $E(F,a,b)$ **кватернионов над полем** $F$[1], если произведение в алгебре $E$ определено согласно правилам

(4.1)

|   | $i$   | $j$  | $k$    |
|---|-------|------|--------|
| $i$ | $a$   | $k$  | $aj$   |
| $j$ | $-k$  | $b$  | $-bi$  |
| $k$ | $-aj$ | $bi$ | $-ab$  |

где $a, b \in F$, $ab \neq 0$. $\square$

Элементы алгебры $E(F,a,b)$ имеют вид

$$x = x^0 + x^1 i + x^2 j + x^3 k$$

где $x^i \in F$, $i = 0, 1, 2, 3$. Кватернион

$$\overline{x} = x^0 - x^1 i - x^2 j - x^3 k$$

называется сопряжённым кватерниону $x$. Мы определим **норму кватерниона** $x$ равенством

(4.2) $$|x|^2 = x\overline{x} = (x^0)^2 - a(x^1)^2 - b(x^2)^2 + ab(x^3)^2$$

Из равенства (4.2) следует, что $E(F,a,b)$ является алгеброй с делением только когда $a < 0$, $b < 0$. Тогда мы можем пронормировать базис так, что $a = -1$, $b = -1$.

Мы будем обозначать символом $E(F)$ алгебру $E(F,-1,-1)$ кватернионов с делением над полем $F$. Мы будем полагать $H = E(R,-1,-1)$. Произведение в алгебре кватернионов $E(F)$ определено согласно правилам

(4.3)

|   | $i$   | $j$  | $k$  |
|---|-------|------|------|
| $i$ | $-1$  | $k$  | $-j$ |
| $j$ | $-k$  | $-1$ | $i$  |
| $k$ | $j$   | $-i$ | $-1$ |

В алгебре $E(F)$ норма кватерниона имеет вид

(4.4) $$|x|^2 = x\overline{x} = (x^0)^2 + (x^1)^2 + (x^2)^2 + (x^3)^2$$

При этом обратный элемент имеет вид

(4.5) $$x^{-1} = |x|^{-2}\overline{x}$$

Внутренний автоморфизм алгебры кватернионов $H$[2]

(4.6) $$p \to qpq^{-1}$$
$$q(ix + jy + kz)q^{-1} = ix' + jy' + kz'$$

---

[1] Я буду следовать определению из [3].
[2] См. [6], с. 643.

описывает вращение вектора с координатами $x$, $y$, $z$. Если $q$ записан в виде суммы скаляра и вектора

$$q = \cos\alpha + (ia + jb + kc)\sin\alpha \quad a^2 + b^2 + c^2 = 1$$

то (4.6) описывает вращение вектора $(x,y,z)$ вокруг вектора $(a,b,c)$ на угол $2\alpha$.

## 5. Башня алгебр

Пусть $F_1$ - алгебра над полем $F_2$. Пусть $\overline{\overline{e}}_{12}$ - базис алгебры $F_1$ над полем $F_2$. Пусть $C_{12\cdot ij}^{\phantom{12\cdot}k}$ - структурные константы алгебры $F_1$ над полем $F_2$.

Пусть $F_2$ - алгебра над полем $F_3$. Пусть $\overline{\overline{e}}_{23}$ - базис алгебры $F_2$ над полем $F_3$. Пусть $C_{23\cdot ij}^{\phantom{23\cdot}k}$ - структурные константы алгебры $F_2$ над полем $F_3$.

Я буду рассматривать алгебру $F_1$ как прямую сумму алгебр $F_2$. Каждое слагаемое суммы я отождествляю с вектором базиса $\overline{\overline{e}}_{12}$. Соответственно, я могу рассматривать алгебру $F_1$ как алгебру над полем $F_3$. Пусть $\overline{\overline{e}}_{13}$ - базис алгебры $F_1$ над полем $F_3$. Индекс базиса $\overline{\overline{e}}_{13}$ состоит из двух индексов: индекса слоя и индекса вектора базиса $\overline{\overline{e}}_{23}$ в слое.

Я буду отождествлять вектор базиса $\overline{e}_{12\cdot i}$ с единицей в соответствующем слое. Тогда

$$(5.1) \qquad \overline{e}_{13\cdot ji} = \overline{e}_{23\cdot j}\overline{e}_{12\cdot i}$$

Произведение векторов базиса $\overline{\overline{e}}_{13}$ имеет вид

$$(5.2) \qquad \overline{e}_{13\cdot ji}\overline{e}_{13\cdot mk} = \overline{e}_{23\cdot j}\overline{e}_{12\cdot i}\overline{e}_{23\cdot m}\overline{e}_{12\cdot k} = C_{23\cdot jm}^{\phantom{23\cdot}a}\overline{e}_{23\cdot a}C_{12\cdot ik}^{\phantom{12\cdot}b}\overline{e}_{12\cdot b}$$

Так как $C_{12\cdot ik}^{\phantom{12\cdot}b} \in F_2$, то разложение $C_{12\cdot ik}^{\phantom{12\cdot}b}$ относительно базиса $\overline{\overline{e}}_{23}$ имеет вид

$$(5.3) \qquad C_{12\cdot ik}^{\phantom{12\cdot}b} = C_{12\cdot ik}^{\phantom{12\cdot}bc}\overline{e}_{23\cdot c}$$

Подставим (5.3) в (5.2)

$$(5.4) \qquad \begin{aligned} \overline{e}_{13\cdot ji}\overline{e}_{13\cdot mk} &= C_{23\cdot jm}^{\phantom{23\cdot}a}\overline{e}_{23\cdot a}C_{12\cdot ik}^{\phantom{12\cdot}bc}\overline{e}_{23\cdot c}\overline{e}_{12\cdot b} \\ &= C_{23\cdot jm}^{\phantom{23\cdot}a}C_{23\cdot ac}^{\phantom{23\cdot}d}C_{12\cdot ik}^{\phantom{12\cdot}bc}\overline{e}_{23\cdot d}\overline{e}_{12\cdot b} \\ &= C_{23\cdot jm}^{\phantom{23\cdot}a}C_{23\cdot ac}^{\phantom{23\cdot}d}C_{12\cdot ik}^{\phantom{12\cdot}bc}\overline{e}_{13\cdot db} \end{aligned}$$

Следовательно, мы можем определить структурные константы алгебры $F_1$ над полем $F_3$

$$(5.5) \qquad C_{13\cdot ji\cdot mk}^{\phantom{13\cdot}\cdot db} = C_{23\cdot jm}^{\phantom{23\cdot}a}C_{23\cdot ac}^{\phantom{23\cdot}d}C_{12\cdot ik}^{\phantom{12\cdot}bc}$$

Чтобы проверить построение, рассмотрим произведение

$$(5.6) \qquad \begin{aligned} \overline{e}_{12\cdot i}\overline{e}_{12\cdot k} &= \overline{e}_{13\cdot 0i}\overline{e}_{13\cdot 0k} \\ &= C_{13\cdot 0i\cdot 0k}^{\phantom{13\cdot}\cdot db}\overline{e}_{13\cdot db} \\ &= C_{23\cdot 00}^{\phantom{23\cdot}a}\,C_{23\cdot ac}^{\phantom{23\cdot}d}C_{12\cdot ik}^{\phantom{12\cdot}bc}\overline{e}_{13\cdot db} \\ &= C_{23\cdot 00}^{\phantom{23\cdot}0}C_{23\cdot 0c}^{\phantom{23\cdot}d}C_{12\cdot ik}^{\phantom{12\cdot}bc}\,\overline{e}_{13\cdot bd} \\ &= C_{23\cdot 0c}^{\phantom{23\cdot}d}C_{12\cdot ik}^{\phantom{12\cdot}bc}\,\overline{e}_{13\cdot bd} \end{aligned}$$

С другой стороны

$$\begin{aligned}
\overline{e}_{12\cdot i}\overline{e}_{12\cdot k} &= C_{12\cdot ik}^{\ \ b}\overline{e}_{12\cdot b} \\
&= C_{12\cdot ik}^{\ \ bd}\overline{e}_{23\cdot d}\overline{e}_{23\cdot 0}\overline{e}_{12\cdot b} \\
&= C_{12\cdot ik}^{\ \ bd}C_{23\cdot 0d}^{\ \ c}\overline{e}_{23\cdot c}\overline{e}_{12\cdot b} \\
&= C_{12\cdot ik}^{\ \ bd}C_{23\cdot 0d}^{\ \ c}\overline{e}_{13\cdot bc}
\end{aligned} \qquad (5.7)$$

Выражения (5.6), (5.7) совпадают.

**Теорема 5.1.** *Если $C_{12\cdot ik}^{\ \ b} \in F_3$, то компоненты $\overline{e}_{12\cdot i}$, $\overline{e}_{23\cdot k}$ вектора $\overline{e}_{13\cdot ik}$ перемножаются независимо*

$$C_{13\cdot\cdot ji\cdot mk}^{\ \ \ \cdot db} = C_{23\cdot jm}^{\ \ d}C_{12\cdot ik}^{\ \ b} \qquad (5.8)$$

*Доказательство.* Из равенства (5.2) следует

$$\overline{e}_{13\cdot ji}\overline{e}_{13\cdot mk} = C_{23\cdot jm}^{\ \ a}\overline{e}_{23\cdot a}C_{12\cdot ik}^{\ \ b}\overline{e}_{12\cdot b} = C_{23\cdot jm}^{\ \ a}C_{12\cdot ik}^{\ \ b}\,\overline{e}_{23\cdot ab} \qquad (5.9)$$

Равенство (5.8) следует из равенства (5.9). $\square$

**Теорема 5.2.** *Пусть $F_1$ - алгебра над полем $F_2$. Пусть $F_2$ - алгебра над полем $F_3$. Координаты $a_{12}^k$ элемента $a \in F_1$ относительно базиса $\overline{\overline{e}}_{12}$ и координаты $a_{13}^{jk}$ элемента $a \in F_1$ относительно базиса $\overline{\overline{e}}_{13}$ связаны равенством*

$$a_{12}^k = a_{13}^{ik}\overline{e}_{23\cdot i} \qquad (5.10)$$

*Доказательство.* Элемент $a \in F_1$ имеет разложение

$$a = a_{12}^k\overline{e}_{12\cdot k} \qquad (5.11)$$

относительно базиса $\overline{\overline{e}}_{12}$ и разложение

$$a = a_{13}^{ik}\overline{e}_{13\cdot ik} \qquad (5.12)$$

относительно базиса $\overline{\overline{e}}_{13}$ Из равенств (5.1), (5.12), следует

$$a = a_{13}^{ik}\overline{e}_{23\cdot i}\overline{e}_{12\cdot k} \qquad (5.13)$$

Равенство (5.10) следует из равенств (5.11),(5.13). $\square$

**Теорема 5.3** (Уравнения Коши-Римана). *Пусть $F_1$ - алгебра над полем $F_2$. Пусть $F_2$ - алгебра над полем $F_3$. Матрица $F_2$-линейного отображения*

$$f : F_1 \to F_1$$

*$F_2$-алгебры $F_1$ удовлетворяет соотношению*

$$f_{jl}^{ik}C_{23\cdot mp}^{\ \ j} = C_{23\cdot mn}^{\ \ i}f_{pl}^{nk} \qquad (5.14)$$

*Доказательство.* Так как $f$ - $F_2$-линейное отображение, то

$$f(ax) = af(x) \quad a \in F_2 \quad x \in F_1 \qquad (5.15)$$

Пусть

$$a = a^j\overline{e}_{23\cdot j} \quad x = x_{13}^{jl}\overline{e}_{13\cdot jl}$$

Тогда

$$f(ax)^{ik} = (ax)^{jl}f_{jl}^{ik} = C_{23\cdot mp}^{\ \ j}a^m x^{pl}f_{jl}^{ik} \qquad (5.16)$$

Аналогично

$$(af(x))^{ik} = C_{23\cdot mn}^{\ \ i}a^m(f(x))^{nk} = C_{23\cdot mn}^{\ \ i}a^m f_{pl}^{nk}x^{pl} \qquad (5.17)$$

Из равенств (5.15), (5.16), (5.17) следует

$$(5.18) \qquad f^{ik}_{jl} C_{23 \cdot mp}^{\phantom{23 \cdot} j} a^m x^{pl} = C_{23 \cdot mn}^{\phantom{23 \cdot} i} a^m f^{nk}_{pl} x^{pl}$$

Так как $a \in F_2$, $x \in F_1$ произвольны, то равенство (5.14) следует из равенства (5.18). □

**Теорема 5.4** (Уравнения Коши-Римана). *Матрица $F_2$-линейного отображения*

$$f : F_2 \to F_2$$

*$F_3$-алгебры $F_2$ удовлетворяет соотношению*

$$(5.19) \qquad f^i_j C_{23 \cdot mp}^{\phantom{23 \cdot} j} = C_{23 \cdot mn}^{\phantom{23 \cdot} i} f^n_p$$

*Доказательство.* Утверждение теоремы будет доказано, если мы в теореме 5.3 положим $F_1 = F_2$. □

**Теорема 5.5** (Уравнения Коши-Римана). *Матрица $C$-линейного отображения*

$$f : C \to C$$

*$R$-алгебры $C$ удовлетворяет соотношению*

$$(5.20) \qquad \begin{aligned} f^0_0 &= \phantom{-}f^1_1 \\ f^0_1 &= -f^1_0 \end{aligned}$$

*Доказательство.* Для поля комплексных чисел $C$ уравнение (5.19) имеет вид

$$(5.21) \qquad \begin{aligned} f^0_0 C_{C \cdot mp}^{\phantom{C \cdot} 0} + f^0_1 C_{C \cdot mp}^{\phantom{C \cdot} 1} &= C_{C \cdot m0}^{\phantom{C \cdot} 0} f^0_p + C_{C \cdot m1}^{\phantom{C \cdot} 0} f^1_p \\ f^1_0 C_{C \cdot mp}^{\phantom{C \cdot} 0} + f^1_1 C_{C \cdot mp}^{\phantom{C \cdot} 1} &= C_{C \cdot m0}^{\phantom{C \cdot} 1} f^0_p + C_{C \cdot m1}^{\phantom{C \cdot} 1} f^1_p \end{aligned}$$

Из уравнений (3.3), (5.21), следует

$$(5.22) \qquad \begin{cases} f^0_0 C_{C \cdot 00}^{\phantom{C \cdot} 0} = C_{C \cdot 00}^{\phantom{C \cdot} 0} f^0_0 \\ f^0_0 C_{C \cdot 11}^{\phantom{C \cdot} 0} = C_{C \cdot 11}^{\phantom{C \cdot} 0} f^1_1 \\ f^0_1 C_{C \cdot 01}^{\phantom{C \cdot} 1} = C_{C \cdot 00}^{\phantom{C \cdot} 0} f^0_1 \\ f^0_1 C_{C \cdot 10}^{\phantom{C \cdot} 1} = C_{C \cdot 11}^{\phantom{C \cdot} 0} f^1_0 \\ f^1_0 C_{C \cdot 00}^{\phantom{C \cdot} 0} = C_{C \cdot 01}^{\phantom{C \cdot} 1} f^1_0 \\ f^1_0 C_{C \cdot 11}^{\phantom{C \cdot} 0} = C_{C \cdot 10}^{\phantom{C \cdot} 1} f^0_1 \\ f^1_1 C_{C \cdot 01}^{\phantom{C \cdot} 1} = C_{C \cdot 01}^{\phantom{C \cdot} 1} f^1_1 \\ f^1_1 C_{C \cdot 10}^{\phantom{C \cdot} 1} = C_{C \cdot 10}^{\phantom{C \cdot} 1} f^0_0 \end{cases}$$

Уравнения (5.20) следует из уравнений (3.3), (5.22). □

Я намеренно записал равенство (5.22) в промежуточной форме, чтобы показать насколько более трудоёмко доказательство теоремы 5.5 по сравнению с доказательством аналогичной теоремы 3.1. Однако мы нигде не пользовались утверждением, что произведение в $F_1$-алгебре $F_2$ коммутативно. Поэтому теорема 5.4 может быть применена к произвольной алгебре.

**Теорема 5.6** (Уравнения Коши-Римана). *Матрица $H$-линейного отображения*

$$f : H \to H$$

$R$-алгебры $H$ удовлетворяет соотношению

$$
\begin{aligned}
f_0^0 &= f_1^1 = f_2^2 = f_3^3 \\
f_1^0 &= -f_0^1 = -f_3^2 = f_2^3 \\
f_2^0 &= f_3^1 = -f_0^2 = -f_1^3 \\
f_3^0 &= -f_2^1 = f_1^2 = -f_0^3
\end{aligned}
\tag{5.23}
$$

*Доказательство.* Согласно теореме [2]-7.3.1. структурные константы алгебры кватернионов имеют вид

$$
\begin{aligned}
C_{H\cdot 00}^{\phantom{H\cdot}0} &= 1 & C_{H\cdot 01}^{\phantom{H\cdot}1} &= 1 & C_{H\cdot 02}^{\phantom{H\cdot}2} &= 1 & C_{H\cdot 03}^{\phantom{H\cdot}3} &= 1 \\
C_{H\cdot 10}^{\phantom{H\cdot}1} &= 1 & C_{H\cdot 11}^{\phantom{H\cdot}0} &= -1 & C_{H\cdot 12}^{\phantom{H\cdot}3} &= 1 & C_{H\cdot 13}^{\phantom{H\cdot}2} &= -1 \\
C_{H\cdot 20}^{\phantom{H\cdot}2} &= 1 & C_{H\cdot 21}^{\phantom{H\cdot}3} &= -1 & C_{H\cdot 22}^{\phantom{H\cdot}0} &= -1 & C_{H\cdot 23}^{\phantom{H\cdot}1} &= 1 \\
C_{H\cdot 30}^{\phantom{H\cdot}3} &= 1 & C_{H\cdot 31}^{\phantom{H\cdot}2} &= 1 & C_{H\cdot 32}^{\phantom{H\cdot}1} &= -1 & C_{H\cdot 33}^{\phantom{H\cdot}0} &= -1
\end{aligned}
\tag{5.24}
$$

Для алгебры кватернионов $H$ уравнение (5.19) имеет вид

$$
\begin{aligned}
&f_0^0 C_{H\cdot mp}^{\phantom{H\cdot}0} + f_1^0 C_{H\cdot mp}^{\phantom{H\cdot}1} + f_2^0 C_{H\cdot mp}^{\phantom{H\cdot}2} + f_3^0 C_{H\cdot mp}^{\phantom{H\cdot}3} \\
&= C_{H\cdot m0}^{\phantom{H\cdot}0} f_p^0 + C_{H\cdot m1}^{\phantom{H\cdot}0} f_p^1 + C_{H\cdot m2}^{\phantom{H\cdot}0} f_p^2 + C_{H\cdot m3}^{\phantom{H\cdot}0} f_p^3 \quad i=0 \\
&f_0^1 C_{H\cdot mp}^{\phantom{H\cdot}0} + f_1^1 C_{H\cdot mp}^{\phantom{H\cdot}1} + f_2^1 C_{H\cdot mp}^{\phantom{H\cdot}2} + f_3^1 C_{H\cdot mp}^{\phantom{H\cdot}3} \\
&= C_{H\cdot m0}^{\phantom{H\cdot}1} f_p^0 + C_{H\cdot m1}^{\phantom{H\cdot}1} f_p^1 + C_{H\cdot m2}^{\phantom{H\cdot}1} f_p^2 + C_{H\cdot m3}^{\phantom{H\cdot}1} f_p^3 \quad i=1 \\
&f_0^2 C_{H\cdot mp}^{\phantom{H\cdot}0} + f_1^2 C_{H\cdot mp}^{\phantom{H\cdot}1} + f_2^2 C_{H\cdot mp}^{\phantom{H\cdot}2} + f_3^2 C_{H\cdot mp}^{\phantom{H\cdot}3} \\
&= C_{H\cdot m0}^{\phantom{H\cdot}2} f_p^0 + C_{H\cdot m1}^{\phantom{H\cdot}2} f_p^1 + C_{H\cdot m2}^{\phantom{H\cdot}2} f_p^2 + C_{H\cdot m3}^{\phantom{H\cdot}2} f_p^3 \quad i=2 \\
&f_0^3 C_{H\cdot mp}^{\phantom{H\cdot}0} + f_1^3 C_{H\cdot mp}^{\phantom{H\cdot}1} + f_2^3 C_{H\cdot mp}^{\phantom{H\cdot}2} + f_3^3 C_{H\cdot mp}^{\phantom{H\cdot}3} \\
&= C_{H\cdot m0}^{\phantom{H\cdot}3} f_p^0 + C_{H\cdot m1}^{\phantom{H\cdot}3} f_p^1 + C_{H\cdot m2}^{\phantom{H\cdot}3} f_p^2 + C_{H\cdot m3}^{\phantom{H\cdot}3} f_p^3 \quad i=3
\end{aligned}
\tag{5.25}
$$

Из уравнений (5.24), (5.25), следует[3]

$$
\begin{cases}
f_1^0 C_{H\cdot 10}^{\phantom{H\cdot}1} = C_{H\cdot 11}^{\phantom{H\cdot}0} f_0^1 & i=0 \ m=1 \ p=0 \\
f_0^0 C_{H\cdot 11}^{\phantom{H\cdot}0} = C_{H\cdot 11}^{\phantom{H\cdot}0} f_1^1 & i=0 \ m=1 \ p=1 \\
f_3^0 C_{H\cdot 12}^{\phantom{H\cdot}3} = C_{H\cdot 11}^{\phantom{H\cdot}0} f_2^1 & i=0 \ m=1 \ p=2 \\
f_2^0 C_{H\cdot 13}^{\phantom{H\cdot}2} = C_{H\cdot 11}^{\phantom{H\cdot}0} f_3^1 & i=0 \ m=1 \ p=3 \\
f_2^0 C_{H\cdot 20}^{\phantom{H\cdot}2} = C_{H\cdot 22}^{\phantom{H\cdot}0} f_0^2 & i=0 \ m=2 \ p=0 \\
f_3^0 C_{H\cdot 21}^{\phantom{H\cdot}3} = C_{H\cdot 22}^{\phantom{H\cdot}0} f_1^2 & i=0 \ m=2 \ p=1 \\
f_0^0 C_{H\cdot 22}^{\phantom{H\cdot}0} = C_{H\cdot 22}^{\phantom{H\cdot}0} f_2^2 & i=0 \ m=2 \ p=2 \\
f_1^0 C_{H\cdot 23}^{\phantom{H\cdot}1} = C_{H\cdot 22}^{\phantom{H\cdot}0} f_3^2 & i=0 \ m=2 \ p=3 \\
f_3^0 C_{H\cdot 30}^{\phantom{H\cdot}3} = C_{H\cdot 33}^{\phantom{H\cdot}0} f_0^3 & i=0 \ m=3 \ p=0 \\
f_2^0 C_{H\cdot 31}^{\phantom{H\cdot}2} = C_{H\cdot 33}^{\phantom{H\cdot}0} f_1^3 & i=0 \ m=3 \ p=1 \\
f_1^0 C_{H\cdot 32}^{\phantom{H\cdot}1} = C_{H\cdot 33}^{\phantom{H\cdot}0} f_2^3 & i=0 \ m=3 \ p=2 \\
f_0^0 C_{H\cdot 33}^{\phantom{H\cdot}0} = C_{H\cdot 33}^{\phantom{H\cdot}0} f_3^3 & i=0 \ m=3 \ p=3
\end{cases}
\tag{5.26}
$$

Уравнения (5.23) следует из уравнений (5.24), (5.26). $\square$

**Теорема 5.7.** *Пусть $F_1$ - алгебра над полем $F_2$. Пусть $F_2$ - алгебра над полем $F_3$. Пусть $F_2$-линейное отображение*

$$f : F_1 \to F_1$$

---

[3] Я не рассматриваю равенства, которые тривиальны или следуют из равенств (5.26).

*имеет представление*

$$f(a)^{\boldsymbol{i}} = f^{\boldsymbol{i}}_{\boldsymbol{j}} a^{\boldsymbol{j}} \quad f^{\boldsymbol{i}}_{\boldsymbol{j}} = f^{\boldsymbol{ik}}_{\boldsymbol{j}} \overline{e}_{23\cdot\boldsymbol{k}} \tag{5.27}$$

*относительно базиса* $\overline{\overline{e}}_{12}$ *и имеет представление*

$$f(a)^{\boldsymbol{ki}} = f^{\boldsymbol{ki}}_{\boldsymbol{lj}} a^{\boldsymbol{lj}} \tag{5.28}$$

*относительно базиса* $\overline{\overline{e}}_{13}$. *Координаты* $f^{\boldsymbol{jk}}_{\boldsymbol{il}}$ *отображения* $f$ *относительно базиса* $\overline{\overline{e}}_{13}$ *имеют вид*

$$f^{\boldsymbol{pi}}_{\boldsymbol{lj}} = C_{23\cdot\boldsymbol{kl}}^{\phantom{23\cdot}\boldsymbol{p}} f^{\boldsymbol{ik}}_{\boldsymbol{j}} \tag{5.29}$$

*Доказательство.* Так как

$$a^{\boldsymbol{j}} = a^{\boldsymbol{lj}} \overline{e}_{23\cdot\boldsymbol{l}}$$

то, из равенства (5.27), следует, что

$$f(a)^{\boldsymbol{i}} = C_{23\cdot\boldsymbol{kl}}^{\phantom{23\cdot}\boldsymbol{p}} f^{\boldsymbol{ik}}_{\boldsymbol{j}} a^{\boldsymbol{lj}} \overline{e}_{23\cdot\boldsymbol{p}} \tag{5.30}$$

Из равенства (5.28), следует

$$f(a)^{\boldsymbol{i}} = f(a)^{\boldsymbol{pi}} \overline{e}_{23\cdot\boldsymbol{p}} = f^{\boldsymbol{pi}}_{\boldsymbol{lj}} a^{\boldsymbol{lj}} \overline{e}_{23\cdot\boldsymbol{p}} \tag{5.31}$$

Из равенств (5.30), (5.31), следует

$$f^{\boldsymbol{pi}}_{\boldsymbol{lj}} a^{\boldsymbol{lj}} \overline{e}_{23\cdot\boldsymbol{p}} = C_{23\cdot\boldsymbol{kl}}^{\phantom{23\cdot}\boldsymbol{p}} f^{\boldsymbol{ik}}_{\boldsymbol{j}} a^{\boldsymbol{lj}} \overline{e}_{23\cdot\boldsymbol{p}} \tag{5.32}$$

Равенство (5.29) следует из равенства (5.32), так как $\overline{\overline{e}}_{23}$ - базис и координаты $a^{\boldsymbol{lj}}$ элемента $a \in F_1$ произвольны. $\square$

## 6. Алгебра кватернионов над полем комплексных чисел

В этом разделе я буду рассматривать алгебру кватернионов $E(C, -1, -1)$, где $C$ - поле комплексных чисел.

Произведение в алгебре $E(C, -1, -1)$ определено согласно таблице

(6.1)

|  | $\overline{e}_{12\cdot\boldsymbol{1}}$ | $\overline{e}_{12\cdot\boldsymbol{2}}$ | $\overline{e}_{12\cdot\boldsymbol{3}}$ |
|---|---|---|---|
| $\overline{e}_{12\cdot\boldsymbol{1}}$ | $-1$ | $\overline{e}_{12\cdot\boldsymbol{3}}$ | $-\overline{e}_{12\cdot\boldsymbol{2}}$ |
| $\overline{e}_{12\cdot\boldsymbol{2}}$ | $-\overline{e}_{12\cdot\boldsymbol{3}}$ | $-1$ | $\overline{e}_{12\cdot\boldsymbol{1}}$ |
| $\overline{e}_{12\cdot\boldsymbol{3}}$ | $\overline{e}_{12\cdot\boldsymbol{2}}$ | $-\overline{e}_{12\cdot\boldsymbol{1}}$ | $-1$ |

Согласно теореме [2]-7.3.1. структурные константы алгебры кватернионов имеют вид

$$\begin{aligned}
&C_{12\cdot\boldsymbol{00}}^{\phantom{12\cdot}\boldsymbol{0}} = 1 \quad C_{12\cdot\boldsymbol{01}}^{\phantom{12\cdot}\boldsymbol{1}} = \phantom{-}1 \quad C_{12\cdot\boldsymbol{02}}^{\phantom{12\cdot}\boldsymbol{2}} = \phantom{-}1 \quad C_{12\cdot\boldsymbol{03}}^{\phantom{12\cdot}\boldsymbol{3}} = \phantom{-}1 \\
&C_{12\cdot\boldsymbol{10}}^{\phantom{12\cdot}\boldsymbol{1}} = 1 \quad C_{12\cdot\boldsymbol{11}}^{\phantom{12\cdot}\boldsymbol{0}} = -1 \quad C_{12\cdot\boldsymbol{12}}^{\phantom{12\cdot}\boldsymbol{3}} = \phantom{-}1 \quad C_{12\cdot\boldsymbol{13}}^{\phantom{12\cdot}\boldsymbol{2}} = -1 \\
&C_{12\cdot\boldsymbol{20}}^{\phantom{12\cdot}\boldsymbol{2}} = 1 \quad C_{12\cdot\boldsymbol{21}}^{\phantom{12\cdot}\boldsymbol{3}} = -1 \quad C_{12\cdot\boldsymbol{22}}^{\phantom{12\cdot}\boldsymbol{0}} = -1 \quad C_{12\cdot\boldsymbol{23}}^{\phantom{12\cdot}\boldsymbol{1}} = \phantom{-}1 \\
&C_{12\cdot\boldsymbol{30}}^{\phantom{12\cdot}\boldsymbol{3}} = 1 \quad C_{12\cdot\boldsymbol{31}}^{\phantom{12\cdot}\boldsymbol{2}} = \phantom{-}1 \quad C_{12\cdot\boldsymbol{32}}^{\phantom{12\cdot}\boldsymbol{1}} = -1 \quad C_{12\cdot\boldsymbol{33}}^{\phantom{12\cdot}\boldsymbol{0}} = -1
\end{aligned}$$

Пусть $\overline{e}_{23\cdot\boldsymbol{k}} = \overline{e}_{C\cdot\boldsymbol{k}}$. Произведение в алгебре $C$ определено согласно правилу (3.2). Согласно теореме 3.1 структурные константы поля комплексных чисел над полем действительных чисел имеют вид (3.3)

Следовательно, алгебра $E(C, -1, -1)$ изоморфна тензорному произведению $C \otimes H$. Поэтому можем выбрать базис

$$\overline{e}_{13\cdot\boldsymbol{ij}} = \overline{e}_{C\cdot\boldsymbol{i}} \otimes \overline{e}_{12\cdot\boldsymbol{j}}$$

(6.2)
$$\overline{e}_{13\cdot\cdot00} = 1 \otimes 1 \quad \overline{e}_{13\cdot\cdot01} = 1 \otimes i \quad \overline{e}_{13\cdot\cdot02} = 1 \otimes j \quad \overline{e}_{13\cdot\cdot03} = 1 \otimes k$$
$$\overline{e}_{13\cdot\cdot10} = i \otimes 1 \quad \overline{e}_{13\cdot\cdot11} = i \otimes i \quad \overline{e}_{13\cdot\cdot12} = i \otimes j \quad \overline{e}_{13\cdot\cdot13} = i \otimes k$$

**Теорема 6.1.** *Таблица произведения алгебры $E(C, -1, -1)$ над полем $R$ имеет вид*

|  | $1 \otimes i$ | $1 \otimes j$ | $1 \otimes k$ | $i \otimes 1$ | $i \otimes i$ | $i \otimes j$ | $i \otimes k$ |
|---|---|---|---|---|---|---|---|
| $1 \otimes i$ | $-1 \otimes 1$ | $1 \otimes k$ | $-1 \otimes j$ | $i \otimes i$ | $-i \otimes 1$ | $i \otimes k$ | $-i \otimes j$ |
| $1 \otimes j$ | $-1 \otimes k$ | $-1 \otimes 1$ | $1 \otimes i$ | $i \otimes j$ | $-i \otimes k$ | $-i \otimes 1$ | $i \otimes i$ |
| $1 \otimes k$ | $1 \otimes j$ | $-1 \otimes i$ | $-1 \otimes 1$ | $i \otimes k$ | $i \otimes j$ | $-i \otimes i$ | $-i \otimes 1$ |
| $i \otimes 1$ | $i \otimes i$ | $i \otimes j$ | $i \otimes k$ | $-1 \otimes 1$ | $-1 \otimes i$ | $-1 \otimes j$ | $-1 \otimes k$ |
| $i \otimes i$ | $-i \otimes 1$ | $i \otimes k$ | $-i \otimes j$ | $-1 \otimes i$ | $1 \otimes 1$ | $-1 \otimes k$ | $1 \otimes j$ |
| $i \otimes j$ | $-i \otimes k$ | $-i \otimes 1$ | $i \otimes i$ | $-1 \otimes j$ | $1 \otimes k$ | $1 \otimes 1$ | $-1 \otimes i$ |
| $i \otimes k$ | $i \otimes j$ | $-i \otimes i$ | $-i \otimes 1$ | $-1 \otimes k$ | $-1 \otimes j$ | $1 \otimes i$ | $1 \otimes 1$ |

*Доказательство.* Таблица записана согласно равенству

$$(\overline{e}_{C\cdot i} \otimes \overline{e}_{12\cdot k})(\overline{e}_{C\cdot j} \otimes \overline{e}_{12\cdot m}) = (\overline{e}_{C\cdot i}\ \overline{e}_{C\cdot j}) \otimes (\overline{e}_{12\cdot k}\ \overline{e}_{12\cdot m})$$

и определению базиса (6.2). □

**Теорема 6.2.** *Структурные константы алгебры $E(C, -1, -1)$ над полем $R$ имеют вид*

$$C_{13\cdot\cdot00\cdot00}^{\cdot00} = 1 \quad C_{13\cdot\cdot00\cdot01}^{\cdot01} = 1 \quad C_{13\cdot\cdot00\cdot02}^{\cdot02} = 1 \quad C_{13\cdot\cdot00\cdot03}^{\cdot03} = 1$$
$$C_{13\cdot\cdot01\cdot00}^{\cdot01} = 1 \quad C_{13\cdot\cdot01\cdot01}^{\cdot00} = -1 \quad C_{13\cdot\cdot01\cdot02}^{\cdot03} = 1 \quad C_{13\cdot\cdot01\cdot03\cdot}^{\cdot02} = -1$$
$$C_{13\cdot\cdot02\cdot00}^{\cdot02} = 1 \quad C_{13\cdot\cdot02\cdot01}^{\cdot03} = -1 \quad C_{13\cdot\cdot02\cdot02}^{\cdot00} = -1 \quad C_{13\cdot\cdot02\cdot03}^{\cdot01} = 1$$
$$C_{13\cdot\cdot03\cdot00}^{\cdot03} = 1 \quad C_{13\cdot\cdot03\cdot01}^{\cdot02} = 1 \quad C_{13\cdot\cdot03\cdot02}^{\cdot01} = -1 \quad C_{13\cdot\cdot03\cdot03}^{\cdot00} = -1$$
$$C_{13\cdot\cdot00\cdot10}^{\cdot10} = 1 \quad C_{13\cdot\cdot00\cdot11}^{\cdot11} = 1 \quad C_{13\cdot\cdot00\cdot12}^{\cdot12} = 1 \quad C_{13\cdot\cdot00\cdot13}^{\cdot13} = 1$$
$$C_{13\cdot\cdot01\cdot10}^{\cdot11} = 1 \quad C_{13\cdot\cdot01\cdot11}^{\cdot10} = -1 \quad C_{13\cdot\cdot01\cdot12}^{\cdot13} = 1 \quad C_{13\cdot\cdot01\cdot13}^{\cdot12} = -1$$
$$C_{13\cdot\cdot02\cdot10}^{\cdot12} = 1 \quad C_{13\cdot\cdot02\cdot11}^{\cdot13} = -1 \quad C_{13\cdot\cdot02\cdot12}^{\cdot10} = -1 \quad C_{13\cdot\cdot02\cdot13}^{\cdot11} = 1$$
$$C_{13\cdot\cdot03\cdot10}^{\cdot13} = 1 \quad C_{13\cdot\cdot03\cdot11}^{\cdot12} = 1 \quad C_{13\cdot\cdot03\cdot12}^{\cdot11} = -1 \quad C_{13\cdot\cdot03\cdot13}^{\cdot10} = -1$$
$$C_{13\cdot\cdot10\cdot00}^{\cdot10} = 1 \quad C_{13\cdot\cdot10\cdot01}^{\cdot11} = 1 \quad C_{13\cdot\cdot10\cdot02}^{\cdot12} = 1 \quad C_{13\cdot\cdot10\cdot03}^{\cdot13} = 1$$
$$C_{13\cdot\cdot11\cdot00}^{\cdot11} = 1 \quad C_{13\cdot\cdot11\cdot01}^{\cdot10} = -1 \quad C_{13\cdot\cdot11\cdot02}^{\cdot13} = 1 \quad C_{13\cdot\cdot11\cdot03}^{\cdot12} = -1$$
$$C_{13\cdot\cdot12\cdot00}^{\cdot12} = 1 \quad C_{13\cdot\cdot12\cdot01}^{\cdot13} = -1 \quad C_{13\cdot\cdot12\cdot02}^{\cdot10} = -1 \quad C_{13\cdot\cdot12\cdot03}^{\cdot11} = 1$$
$$C_{13\cdot\cdot13\cdot00}^{\cdot13} = 1 \quad C_{13\cdot\cdot13\cdot01}^{\cdot12} = 1 \quad C_{13\cdot\cdot13\cdot02}^{\cdot11} = -1 \quad C_{13\cdot\cdot13\cdot03}^{\cdot10} = -1$$
$$C_{13\cdot\cdot10\cdot10}^{\cdot00} = -1 \quad C_{13\cdot\cdot10\cdot11}^{\cdot01} = -1 \quad C_{13\cdot\cdot10\cdot12}^{\cdot02} = -1 \quad C_{13\cdot\cdot10\cdot13}^{\cdot03} = -1$$
$$C_{13\cdot\cdot11\cdot10}^{\cdot01} = -1 \quad C_{13\cdot\cdot11\cdot11}^{\cdot00} = 1 \quad C_{13\cdot\cdot11\cdot12}^{\cdot03} = -1 \quad C_{13\cdot\cdot11\cdot13}^{\cdot02} = 1$$
$$C_{13\cdot\cdot12\cdot10}^{\cdot02} = -1 \quad C_{13\cdot\cdot12\cdot11}^{\cdot03} = 1 \quad C_{13\cdot\cdot12\cdot12}^{\cdot00} = 1 \quad C_{13\cdot\cdot12\cdot13}^{\cdot01} = -1$$
$$C_{13\cdot\cdot13\cdot10}^{\cdot03} = -1 \quad C_{13\cdot\cdot13\cdot11}^{\cdot02} = -1 \quad C_{13\cdot\cdot13\cdot12}^{\cdot01} = 1 \quad C_{13\cdot\cdot13\cdot13}^{\cdot00} = 1$$

*Доказательство.* Утверждение теоремы можно рассматривать либо как следствие теоремы 6.1, либо как следствие теоремы 5.1. □

**Теорема 6.3** (Уравнения Коши-Римана). *Матрица линейной функции*
$$y^{ik} = x^{jm} f^{ik}_{jm}$$
*алгебры $E(C,-1,-1)$ удовлетворяет соотношению*

(6.3)
$$\begin{aligned} f^{\cdot 0j}_{\cdot 0i} &= f^{\cdot 1j}_{\cdot 1i} \\ f^{\cdot 1j}_{\cdot 0i} &= -f^{\cdot 0j}_{\cdot 1i} \end{aligned}$$

*Доказательство.* Из равенств [2]-(3.1.17), (5.8), (3.3) следует

(6.4)
$$\begin{aligned} f^{\cdot 0j}_{\cdot 0i} &= f^{\cdot ka \cdot rc} \, C_{13 \cdot \cdot ka \cdot 0i}^{\cdot pb} \, C_{13 \cdot \cdot pb \cdot rc}^{\cdot 0j} \\ &= f^{\cdot ka \cdot rc} \, C_{C \cdot k0}^{\ p} \, C_{12 \cdot ai}^{\ b} \, C_{C \cdot pr}^{\ 0} \, C_{12 \cdot bc}^{\ j} \\ &= f^{\cdot 0a \cdot 0c} \, C_{C \cdot 00}^{\ 0} \, C_{12 \cdot ai}^{\ b} \, C_{C \cdot 00}^{\ 0} \, C_{12 \cdot bc}^{\ j} \\ &\quad + f^{\cdot 1a \cdot 1c} \, C_{C \cdot 10}^{\ 1} \, C_{12 \cdot ai}^{\ b} \, C_{C \cdot 11}^{\ 0} \, C_{12 \cdot bc}^{\ j} \\ &= f^{\cdot 0a \cdot 0c} \, C_{12 \cdot ai}^{\ b} \, C_{12 \cdot bc}^{\ j} - f^{\cdot 1a \cdot 1c} \, C_{12 \cdot ai}^{\ b} \, C_{12 \cdot bc}^{\ j} \\ &= (f^{\cdot 0a \cdot 0c} - f^{\cdot 1a \cdot 1c}) \, C_{12 \cdot ai}^{\ b} \, C_{12 \cdot bc}^{\ j} \end{aligned}$$

(6.5)
$$\begin{aligned} f^{\cdot 1j}_{\cdot 0i} &= f^{\cdot ka \cdot rc} \, C_{13 \cdot \cdot ka \cdot 0i}^{\cdot pb} \, C_{13 \cdot \cdot pb \cdot rc}^{\cdot 1j} \\ &= f^{\cdot ka \cdot rc} \, C_{C \cdot k0}^{\ p} \, C_{12 \cdot ai}^{\ b} \, C_{C \cdot pr}^{\ 1} \, C_{12 \cdot bc}^{\ j} \\ &= f^{\cdot 0a \cdot 1c} \, C_{C \cdot 00}^{\ 0} \, C_{12 \cdot ai}^{\ b} \, C_{C \cdot 01}^{\ 1} \, C_{12 \cdot bc}^{\ j} \\ &\quad + f^{\cdot 1a \cdot 0c} \, C_{C \cdot 10}^{\ 1} \, C_{12 \cdot ai}^{\ b} \, C_{C \cdot 10}^{\ 1} \, C_{12 \cdot bc}^{\ j} \\ &= f^{\cdot 0a \cdot 1c} \, C_{12 \cdot ai}^{\ b} \, C_{12 \cdot bc}^{\ j} + f^{\cdot 1a \cdot 0c} \, C_{12 \cdot ai}^{\ b} \, C_{12 \cdot bc}^{\ j} \\ &= (f^{\cdot 0a \cdot 1c} + f^{\cdot 1a \cdot 0c}) \, C_{12 \cdot ai}^{\ b} \, C_{12 \cdot bc}^{\ j} \end{aligned}$$

(6.6)
$$\begin{aligned} f^{\cdot 0j}_{\cdot 1i} &= f^{\cdot ka \cdot rc} \, C_{13 \cdot \cdot ka \cdot 1i}^{\cdot pb} \, C_{13 \cdot \cdot pb \cdot rc}^{\cdot 0j} \\ &= f^{\cdot ka \cdot rc} \, C_{C \cdot k1}^{\ p} \, C_{12 \cdot ai}^{\ b} \, C_{C \cdot pr}^{\ 0} \, C_{12 \cdot bc}^{\ j} \\ &= f^{\cdot 0a \cdot 1c} \, C_{C \cdot 01}^{\ 1} \, C_{12 \cdot ai}^{\ b} \, C_{C \cdot 11}^{\ 0} \, C_{12 \cdot bc}^{\ j} \\ &\quad + f^{\cdot 1a \cdot 0c} \, C_{C \cdot 11}^{\ 0} \, C_{12 \cdot ai}^{\ b} \, C_{C \cdot 00}^{\ 0} \, C_{12 \cdot bc}^{\ j} \\ &= -f^{\cdot 0a \cdot 1c} \, C_{12 \cdot ai}^{\ b} \, C_{12 \cdot bc}^{\ j} - f^{\cdot 1a \cdot 0c} \, C_{12 \cdot ai}^{\ b} \, C_{12 \cdot bc}^{\ j} \\ &= -(f^{\cdot 0a \cdot 1c} + f^{\cdot 1a \cdot 0c}) \, C_{12 \cdot ai}^{\ b} \, C_{12 \cdot bc}^{\ j} \end{aligned}$$

(6.7)
$$\begin{aligned} f^{\cdot 1j}_{\cdot 1i} &= f^{\cdot ka \cdot rc} \, C_{13 \cdot \cdot ka \cdot 1i}^{\cdot pb} \, C_{13 \cdot \cdot pb \cdot rc}^{\cdot 1j} \\ &= f^{\cdot ka \cdot rc} \, C_{C \cdot k1}^{\ p} \, C_{12 \cdot ai}^{\ b} \, C_{C \cdot pr}^{\ 1} \, C_{12 \cdot bc}^{\ j} \\ &= f^{\cdot 0a \cdot 0c} \, C_{C \cdot 01}^{\ 1} \, C_{12 \cdot ai}^{\ b} \, C_{C \cdot 10}^{\ 1} \, C_{12 \cdot bc}^{\ j} \\ &\quad + f^{\cdot 1a \cdot 1c} \, C_{C \cdot 11}^{\ 0} \, C_{12 \cdot ai}^{\ b} \, C_{C \cdot 01}^{\ 1} \, C_{12 \cdot bc}^{\ j} \\ &= f^{\cdot 0a \cdot 0c} \, C_{12 \cdot ai}^{\ b} \, C_{12 \cdot bc}^{\ j} - f^{\cdot 1a \cdot 1c} \, C_{12 \cdot ai}^{\ b} \, C_{12 \cdot bc}^{\ j} \\ &= (f^{\cdot 0a \cdot 0c} - f^{\cdot 1a \cdot 1c}) \, C_{12 \cdot ai}^{\ b} \, C_{12 \cdot bc}^{\ j} \end{aligned}$$

Равенство (6.3) следует из сравнения равенств (6.4) и (6.7), (6.5) и (6.6). $\square$

**Теорема 6.4** (Уравнения Коши-Римана). *Если матрица*

$$\left(\frac{\partial y^{\cdot ij}}{\partial x^{\cdot kl}}\right)$$

*является матрицей Якоби функции в алгебре* $E(C,-1,-1)$, *то*

$$\frac{\partial y^{\cdot 1j}}{\partial x^{\cdot 0i}} = -\frac{\partial y^{\cdot 0j}}{\partial x^{\cdot 1i}}$$
$$\frac{\partial y^{\cdot 0j}}{\partial x^{\cdot 0i}} = \frac{\partial y^{\cdot 1j}}{\partial x^{\cdot 1i}}$$

*Доказательство.* Следствие теоремы 6.3. □

## 7. Алгебра $C \otimes (C \otimes H)$

Алгебра $C \otimes (C \otimes H)$ не является алгеброй кватернионов. Однако я рассматриваю здесь эту алгебру так как в некотором смысле эта алгебра похожа на алгебру $E(C,-1,-1) = C \otimes H$.

Алгебра $C \otimes (C \otimes H)$ интересна и с другой точки зрения. При рассмотрении этой алгебры становится очевидным, что тензорное произведение некоммутативных колец неассоциативно. Так как $C$ - поле, то $C \otimes C$ изоморфно $C$, и следовательно, $(C \otimes C) \otimes H$ изоморфно $C \otimes H$. Однако, так как $H$ не является алгеброй над полем комплексных чисел, то алгебра $C \otimes (C \otimes H)$ отличается от алгебры $(C \otimes C) \otimes H = C \otimes H$.

Алгебра $C \otimes (C \otimes H)$ имеет настолько большую размерность, что становится нецелесообразным рассматривать таблицу умножения и структурные константы этой алгебры. Однако нетрудно увидеть, что структура этой алгебры во многом похожа на структуру алгебр, рассмотренных в разделе 5.

Мы будем представлять базис алгебры $C \otimes (C \otimes H)$ в виде

(7.1) $$\overline{\overline{e}}_{kji} = \overline{e}_{C \cdot k} \otimes (\overline{e}_{C \cdot j} \otimes \overline{e}_{H \cdot i})$$

где $\overline{\overline{e}}_C$ - базис алгебры $C$ и $\overline{\overline{e}}_H$ - базис алгебры $H$. Соответственно, произведение в алгебре $C \otimes (C \otimes H)$ определено покомпонентно

(7.2) $$a_1 \otimes (a_2 \otimes a_3) \; b_1 \otimes (b_2 \otimes b_2) = (a_1 b_1) \otimes ((a_2 b_2) \otimes (a_3 b_3))$$

**Теорема 7.1.** *Структурные константы алгебры* $C \otimes (C \otimes H)$ *имеют вид*

$$B_{\cdot pji \cdot qmk}^{\cdot rdb} = C_{C \cdot pq}^{\cdot r} \; C_{C \cdot jm}^{\cdot d} \; C_{H \cdot ik}^{\cdot b}$$

*где* $C_{C \cdot pq}^{\cdot r}$ - *структурные константы поля комплексных чисел,* $C_{H \cdot ik}^{\cdot b}$ - *структурные константы алгебры кватернионов.*

*Доказательство.* Для доказательства теоремы достаточно сравнить следующие равенства

$$\begin{aligned}
\overline{e}_{pji} \; \overline{e}_{qmk} &= B_{\cdot pji \cdot qmk}^{\cdot rdb} \; \overline{e}_{rdb}\\
\overline{e}_{pji} \; \overline{e}_{qmk} &= \overline{e}_{C \cdot p} \otimes (\overline{e}_{C \cdot j} \otimes \overline{e}_{H \cdot i}) \overline{e}_{C \cdot q} \otimes (\overline{e}_{C \cdot m} \otimes \overline{e}_{H \cdot k}))\\
&= (\overline{e}_{C \cdot p} \; \overline{e}_{C \cdot q}) \otimes ((\overline{e}_{C \cdot j} \; \overline{e}_{C \cdot m}) \otimes (\overline{e}_{H \cdot i}) \; \overline{e}_{H \cdot k}))\\
&= C_{C \cdot pq}^{\cdot r} C_{C \cdot jm}^{\cdot d} C_{H \cdot ik}^{\cdot b} \overline{e}_{C \cdot r} \otimes (\overline{e}_{C \cdot d} \otimes \overline{e}_{H \cdot b})\\
&= C_{C \cdot pq}^{\cdot r} \; C_{C \cdot jm}^{\cdot d} \; C_{H \cdot ik}^{\cdot b} \; \overline{e}_{rdb}
\end{aligned}$$

□

**Теорема 7.2.** *Структурные константы алгебры* $C \otimes (C \otimes H)$ *имеют вид*

$$
\begin{array}{l|l}
B^{\cdot 00b}_{\cdot 00i \cdot 00k} = \phantom{-}C_{H \cdot ik}^{\phantom{H} \cdot b} & B^{\cdot 01b}_{\cdot 00i \cdot 01k} = \phantom{-}C_{H \cdot ik}^{\phantom{H} \cdot b} \\
B^{\cdot 01b}_{\cdot 01i \cdot 00k} = \phantom{-}C_{H \cdot ik}^{\phantom{H} \cdot b} & B^{\cdot 01b}_{\cdot 01i \cdot 01k} = -C_{H \cdot ik}^{\phantom{H} \cdot b} \\
B^{\cdot 10b}_{\cdot 00i \cdot 10k} = \phantom{-}C_{H \cdot ik}^{\phantom{H} \cdot b} & B^{\cdot 11b}_{\cdot 00i \cdot 11k} = \phantom{-}C_{H \cdot ik}^{\phantom{H} \cdot b} \\
B^{\cdot 11b}_{\cdot 01i \cdot 10k} = \phantom{-}C_{H \cdot ik}^{\phantom{H} \cdot b} & B^{\cdot 11b}_{\cdot 01i \cdot 11k} = -C_{H \cdot ik}^{\phantom{H} \cdot b} \\
B^{\cdot 10b}_{\cdot 10i \cdot 00k} = \phantom{-}C_{H \cdot ik}^{\phantom{H} \cdot b} & B^{\cdot 11b}_{\cdot 10i \cdot 01k} = \phantom{-}C_{H \cdot ik}^{\phantom{H} \cdot b} \\
B^{\cdot 11b}_{\cdot 11i \cdot 00k} = \phantom{-}C_{H \cdot ik}^{\phantom{H} \cdot b} & B^{\cdot 11b}_{\cdot 11i \cdot 01k} = -C_{H \cdot ik}^{\phantom{H} \cdot b} \\
B^{\cdot 00b}_{\cdot 10i \cdot 10k} = -C_{H \cdot ik}^{\phantom{H} \cdot b} & B^{\cdot 01b}_{\cdot 10i \cdot 11k} = -C_{H \cdot ik}^{\phantom{H} \cdot b} \\
B^{\cdot 01b}_{\cdot 11i \cdot 10k} = -C_{H \cdot ik}^{\phantom{H} \cdot b} & B^{\cdot 01b}_{\cdot 11i \cdot 11k} = \phantom{-}C_{H \cdot ik}^{\phantom{H} \cdot b}
\end{array}
\tag{7.3}
$$

*Доказательство.* Утверждение теоремы является следствием теорем 7.1, 3.1. □

**Теорема 7.3** (Уравнения Коши-Римана)**.** *Матрица линейной функции*

$$y^{ikp} = x^{jmr} f^{ikp}_{jmr}$$

*алгебры* $C \otimes (C \otimes H)$ *удовлетворяет соотношению*

$$
\begin{array}{llll}
f^{\cdot 00j}_{\cdot 00i} = & f^{\cdot 10j}_{\cdot 10i} = & f^{\cdot 01j}_{\cdot 01i} = & f^{\cdot 11j}_{\cdot 11i} \\
f^{\cdot 01j}_{\cdot 00i} = & f^{\cdot 11j}_{\cdot 10i} = -f^{\cdot 00j}_{\cdot 01i} = -f^{\cdot 10j}_{\cdot 11i} \\
f^{\cdot 10j}_{\cdot 00i} = -f^{\cdot 00j}_{\cdot 10i} = & f^{\cdot 11j}_{\cdot 01i} = -f^{\cdot 01j}_{\cdot 11i} \\
f^{\cdot 11j}_{\cdot 00i} = -f^{\cdot 01j}_{\cdot 10i} = -f^{\cdot 10j}_{\cdot 01i} = & f^{\cdot 00j}_{\cdot 11i}
\end{array}
\tag{7.4}
$$

*Доказательство.* Опираясь на равенства [2]-(3.1.17), (7.3), мы можем повторить вычисления, выполненные для доказательства теоремы 6.3. Однако очевидно, что эти вычисления отличаются только набором индексов. Таким образом, из теоремы 6.3 следует

$$
\begin{aligned}
f^{\cdot 0mk}_{\cdot 0ji} &= \phantom{-}f^{\cdot 1mk}_{\cdot 1ji} \\
f^{\cdot 1mk}_{\cdot 0ji} &= -f^{\cdot 0mk}_{\cdot 1ji}
\end{aligned}
\tag{7.5}
$$

$$
\begin{aligned}
f^{\cdot q0k}_{\cdot p0i} &= \phantom{-}f^{\cdot q1k}_{\cdot p1i} \\
f^{\cdot q1k}_{\cdot p0i} &= -f^{\cdot q0k}_{\cdot p1i}
\end{aligned}
\tag{7.6}
$$

Равенства (7.4) следуют из равенств (7.5), (7.6). □

**Теорема 7.4** (Уравнения Коши-Римана)**.** *Если матрица*

$$\left( \frac{\partial y^{\cdot ijp}}{\partial x^{\cdot klq}} \right)$$

*является матрицей Якоби функции в алгебре* $C \otimes (C \otimes H)$, *то*

$$\frac{\partial y^{\cdot 00i}}{\partial x^{\cdot 00j}} = \frac{\partial y^{\cdot 10i}}{\partial x^{\cdot 10j}} = \frac{\partial y^{\cdot 01i}}{\partial x^{\cdot 01j}} = \frac{\partial y^{\cdot 11i}}{\partial x^{\cdot 11j}}$$

$$\frac{\partial y^{\cdot 00i}}{\partial x^{\cdot 01j}} = \frac{\partial y^{\cdot 10i}}{\partial x^{\cdot 11j}} = -\frac{\partial y^{\cdot 01i}}{\partial x^{\cdot 00j}} = -\frac{\partial y^{\cdot 11i}}{\partial x^{\cdot 10j}}$$

$$\frac{\partial y^{\cdot 00i}}{\partial x^{\cdot 10j}} = -\frac{\partial y^{\cdot 10i}}{\partial x^{\cdot 00j}} = \frac{\partial y^{\cdot 01i}}{\partial x^{\cdot 11j}} = -\frac{\partial y^{\cdot 11i}}{\partial x^{\cdot 01j}}$$

$$\frac{\partial y^{\cdot 00i}}{\partial x^{\cdot 11j}} = -\frac{\partial y^{\cdot 10i}}{\partial x^{\cdot 01j}} = -\frac{\partial y^{\cdot 01i}}{\partial x^{\cdot 10j}} = \frac{\partial y^{\cdot 11i}}{\partial x^{\cdot 00j}}$$

*Доказательство.* Следствие теоремы 7.3. □

## 8. Алгебра кватернионов $E(R, a, b)$

Положим $\overline{e}_0 = 1$, $\overline{e}_1 = i$, $\overline{e}_2 = j$, $\overline{e}_3 = k$. Согласно равенству (4.1) структурные константы алгебры $E(R, a, b)$ имеют вид

$$\begin{aligned}
C_{00}^0 &= 1 & C_{01}^1 &= 1 & C_{02}^2 &= 1 & C_{03}^3 &= 1 \\
C_{10}^1 &= 1 & C_{11}^0 &= a & C_{12}^3 &= 1 & C_{13}^2 &= a \\
C_{20}^2 &= 1 & C_{21}^3 &= -1 & C_{22}^0 &= b & C_{23}^1 &= -b \\
C_{30}^3 &= 1 & C_{31}^2 &= -a & C_{32}^1 &= b & C_{33}^0 &= -ab
\end{aligned}$$

**Теорема 8.1.** *Стандартные компоненты линейной функции и координаты соответствующего линейного преобразования над полем $R$ удовлетворяют соотношениям*

(8.1)
$$\begin{cases}
f_0^0 = f^{00} + af^{11} + bf^{22} - abf^{33} \\
f_1^1 = f^{00} + af^{11} - bf^{22} + abf^{33} \\
f_2^2 = f^{00} - af^{11} + bf^{22} + abf^{33} \\
f_3^3 = f^{00} - af^{11} - bf^{22} - abf^{33}
\end{cases}$$

(8.2)
$$\begin{cases}
f_0^1 = f^{01} + f^{10} - bf^{23} + bf^{32} \\
f_1^0 = af^{01} + af^{10} + abf^{23} - abf^{32} \\
f_2^3 = -f^{01} + f^{10} + bf^{23} + bf^{32} \\
f_3^2 = -af^{01} + af^{10} - abf^{23} - abf^{32}
\end{cases}$$

(8.3)
$$\begin{cases}
f_0^2 = f^{02} + af^{13} + f^{20} - af^{31} \\
f_1^3 = f^{02} + af^{13} - f^{20} + af^{31} \\
f_2^0 = bf^{02} - abf^{13} + bf^{20} + abf^{31} \\
f_3^1 = bf^{02} - abf^{13} - bf^{20} - abf^{31}
\end{cases}$$

$$(8.4)\begin{cases} f^3_0 = \phantom{-a}f^{03}+\phantom{a}f^{12}-\phantom{a}f^{21}+\phantom{a}f^{30} \\ f^2_1 = \phantom{-}af^{03}+af^{12}+af^{21}-af^{30} \\ f^1_2 =-\phantom{a}bf^{03}+bf^{12}+bf^{21}+bf^{30} \\ f^0_3 =-abf^{03}+abf^{12}-abf^{21}-abf^{30} \end{cases}$$

$$(8.5)\begin{cases} f^{00} = \phantom{-}\dfrac{1}{4}\phantom{a}f^0_0+\dfrac{1}{4}\phantom{a}f^1_1+\dfrac{1}{4}\phantom{a}f^2_2+\dfrac{1}{4}\phantom{a}f^3_3 \\ f^{11} = \phantom{-}\dfrac{1}{4a}\phantom{a}f^0_0+\dfrac{1}{4a}f^1_1-\dfrac{1}{4a}f^2_2-\dfrac{1}{4a}f^3_3 \\ f^{22} = \phantom{-}\dfrac{1}{4b}\phantom{a}f^0_0-\dfrac{1}{4b}f^1_1+\dfrac{1}{4b}f^2_2-\dfrac{1}{4b}f^3_3 \\ f^{33} =-\dfrac{1}{4ab}f^0_0+\dfrac{1}{4ab}f^1_1+\dfrac{1}{4ab}f^2_2-\dfrac{1}{4ab}f^3_3 \end{cases}$$

$$(8.6)\begin{cases} f^{10} = \phantom{-}\dfrac{1}{4a}f^0_1+\dfrac{1}{4}f^1_0+\dfrac{1}{4a}f^2_3+\dfrac{1}{4}f^3_2 \\ f^{01} = \phantom{-}\dfrac{1}{4a}f^0_1+\dfrac{1}{4}f^1_0-\dfrac{1}{4a}f^2_3-\dfrac{1}{4}f^3_2 \\ f^{32} =-\dfrac{1}{4ab}f^0_1+\dfrac{1}{4b}f^1_0-\dfrac{1}{4ab}f^2_3+\dfrac{1}{4b}f^3_2 \\ f^{23} = \phantom{-}\dfrac{1}{4ab}f^0_1-\dfrac{1}{4b}f^1_0-\dfrac{1}{4ab}f^2_3+\dfrac{1}{4b}f^3_2 \end{cases}$$

$$(8.7)\begin{cases} f^{20} = \phantom{-}\dfrac{1}{4b}f^0_2-\dfrac{1}{4b}f^1_3+\dfrac{1}{4}f^2_0-\dfrac{1}{4}f^3_1 \\ f^{31} = \phantom{-}\dfrac{1}{4ab}f^0_2-\dfrac{1}{4ab}f^1_3-\dfrac{1}{4a}f^2_0+\dfrac{1}{4a}f^3_1 \\ f^{02} = \phantom{-}\dfrac{1}{4b}f^0_2+\dfrac{1}{4b}f^1_3+\dfrac{1}{4}f^2_0+\dfrac{1}{4}f^3_1 \\ f^{13} =-\dfrac{1}{4ab}f^0_2-\dfrac{1}{4ab}f^1_3+\dfrac{1}{4a}f^2_0+\dfrac{1}{4a}f^3_1 \end{cases}$$

$$(8.8)\begin{cases} f^{30} =-\dfrac{1}{4ab}f^0_3+\dfrac{1}{4b}f^1_2-\dfrac{1}{4a}f^2_1+\dfrac{1}{4}f^3_0 \\ f^{21} =-\dfrac{1}{4ab}f^0_3+\dfrac{1}{4b}f^1_2+\dfrac{1}{4a}f^2_1-\dfrac{1}{4}f^3_0 \\ f^{12} = \phantom{-}\dfrac{1}{4ab}f^0_3+\dfrac{1}{4b}f^1_2+\dfrac{1}{4a}f^2_1+\dfrac{1}{4}f^3_0 \\ f^{03} =-\dfrac{1}{4ab}f^0_3-\dfrac{1}{4b}f^1_2+\dfrac{1}{4a}f^2_1+\dfrac{1}{4}f^3_0 \end{cases}$$

*Доказательство.* Пользуясь равенством [2]-3.1.17 получаем соотношения

$$\begin{aligned} f^0_0 &= f^{kr}\ C^p_{k0}\ C^0_{pr} \\ &= f^{00}\ C^0_{00}\ C^0_{00} + f^{11}\ C^1_{10}\ C^0_{11} + f^{22}\ C^2_{20}\ C^0_{22} + f^{33}\ C^3_{30}\ C^0_{33} \\ &= f^{00} + af^{11} + bf^{22} - abf^{33} \end{aligned}$$

$$\begin{aligned} f^1_0 &= f^{kr}\ C^p_{k0}\ C^1_{pr} \\ &= f^{01}\ C^0_{00}\ C^1_{01} + f^{10}\ C^1_{10}\ C^1_{10} + f^{23}\ C^2_{20}\ C^1_{23} + f^{32}\ C^3_{30}\ C^1_{32} \\ &= f^{01} + f^{10} - bf^{23} + bf^{32} \end{aligned}$$

$$f_0^2 = f^{kr}\ C_{k0}^p\ C_{pr}^2$$
$$= f^{02}\ C_{00}^0\ C_{02}^2 + f^{13}\ C_{10}^1\ C_{13}^2 + f^{20}\ C_{20}^2\ C_{20}^2 + f^{31}\ C_{30}^3\ C_{31}^2$$
$$= f^{02} + af^{13} + f^{20} - af^{31}$$

$$f_0^3 = f^{kr}\ C_{k0}^p\ C_{pr}^3$$
$$= f^{03}\ C_{00}^0\ C_{03}^3 + f^{12}\ C_{10}^1\ C_{12}^3 + f^{21}\ C_{20}^2\ C_{21}^3 + f^{30}\ C_{30}^3\ C_{30}^3$$
$$= f^{03} + f^{12} - f^{21} + f^{30}$$

$$f_1^0 = f^{kr}\ C_{k1}^p\ C_{pr}^0$$
$$= f^{01}\ C_{01}^1\ C_{11}^0 + f^{10}\ C_{11}^0\ C_{00}^0 + f^{23}\ C_{21}^3\ C_{33}^0 + f^{32}\ C_{31}^2\ C_{22}^0$$
$$= af^{01} + af^{10} + abf^{23} - abf^{32}$$

$$f_1^1 = f^{kr}\ C_{k1}^p\ C_{pr}^1$$
$$= f^{00}\ C_{01}^1\ C_{10}^1 + f^{11}\ C_{11}^0\ C_{01}^1 + f^{22}\ C_{21}^3\ C_{32}^1 + f^{33}\ C_{31}^2\ C_{23}^1$$
$$= f^{00} + af^{11} - bf^{22} + abf^{33}$$

$$f_1^2 = f^{kr}\ C_{k1}^p\ C_{pr}^2$$
$$= f^{03}\ C_{01}^1\ C_{13}^2 + f^{12}\ C_{11}^0\ C_{02}^2 + f^{21}\ C_{21}^3\ C_{31}^2 + f^{30}\ C_{31}^2\ C_{20}^2$$
$$= af^{03} + af^{12} + af^{21} - af^{30}$$

$$f_1^3 = f^{kr}\ C_{k1}^p\ C_{pr}^3$$
$$= f^{02}\ C_{01}^1\ C_{12}^3 + f^{13}\ C_{11}^0\ C_{03}^3 + f^{20}\ C_{21}^3\ C_{30}^3 + f^{31}\ C_{31}^2\ C_{21}^3$$
$$= f^{02} + af^{13} - f^{20} + af^{31}$$

$$f_2^0 = f^{kr}\ C_{k2}^p\ C_{pr}^0$$
$$= f^{02}\ C_{02}^2\ C_{22}^0 + f^{13}\ C_{12}^3\ C_{33}^0 + f^{20}\ C_{22}^0\ C_{00}^0 + f^{31}\ C_{32}^1\ C_{11}^0$$
$$= bf^{02} - abf^{13} + bf^{20} + abf^{31}$$

$$f_2^1 = f^{kr}\ C_{k2}^p\ C_{pr}^1$$
$$= f^{03}\ C_{02}^2\ C_{23}^1 + f^{12}\ C_{12}^3\ C_{32}^1 + f^{21}\ C_{22}^0\ C_{01}^1 + f^{30}\ C_{32}^1\ C_{10}^1$$
$$= -bf^{03} + bf^{12} + bf^{21} + bf^{30}$$

$$f_2^2 = f^{kr}\ C_{k2}^p\ C_{pr}^2$$
$$= f^{00}\ C_{02}^2\ C_{20}^2 + f^{11}\ C_{12}^3\ C_{31}^2 + f^{22}\ C_{22}^0\ C_{02}^2 + f^{33}\ C_{32}^1\ C_{13}^2$$
$$= f^{00} - af^{11} + bf^{22} + abf^{33}$$

$$f_2^3 = f^{kr}\ C_{k2}^p\ C_{pr}^3$$
$$= f^{01}\ C_{02}^2\ C_{21}^3 + f^{10}\ C_{12}^3\ C_{30}^3 + f^{23}\ C_{22}^0\ C_{03}^3 + f^{32}\ C_{32}^1\ C_{12}^3$$
$$= -f^{01} + f^{10} + bf^{23} + bf^{32}$$

$$\begin{aligned}
f_3^0 &= f^{kr}\ C_{k3}^p\ C_{pr}^0 \\
&= f^{03}\ C_{03}^3\ C_{33}^0 + f^{12}\ C_{13}^2\ C_{22}^0 + f^{21}\ C_{23}^1\ C_{11}^0 + f^{30}\ C_{33}^0\ C_{00}^0 \\
&= -abf^{03} + abf^{12} - abf^{21} - abf^{30}
\end{aligned}$$

$$\begin{aligned}
f_3^1 &= f^{kr}\ C_{k3}^p\ C_{pr}^1 \\
&= f^{02}\ C_{03}^3\ C_{32}^1 + f^{13}\ C_{13}^2\ C_{23}^1 + f^{20}\ C_{23}^1\ C_{10}^1 + f^{31}\ C_{33}^0\ C_{01}^1 \\
&= bf^{02} - abf^{13} - bf^{20} - abf^{31}
\end{aligned}$$

$$\begin{aligned}
f_3^2 &= f^{kr}\ C_{k3}^p\ C_{pr}^2 \\
&= f^{01}\ C_{03}^3\ C_{31}^2 + f^{10}\ C_{13}^2\ C_{20}^2 + f^{23}\ C_{23}^1\ C_{13}^2 + f^{32}\ C_{33}^0\ C_{02}^2 \\
&= -af^{01} + af^{10} - abf^{23} - abf^{32}
\end{aligned}$$

$$\begin{aligned}
f_3^3 &= f^{kr}\ C_{k3}^p\ C_{pr}^3 \\
&= f^{00}\ C_{03}^3\ C_{30}^3 + f^{11}\ C_{13}^2\ C_{21}^3 + f^{22}\ C_{23}^1\ C_{12}^3 + f^{33}\ C_{33}^0\ C_{03}^3 \\
&= f^{00} - af^{11} - bf^{22} - abf^{33}
\end{aligned}$$

Мы группируем эти соотношения в системы линейных уравнений (8.1), (8.2), (8.3), (8.4).

(8.5) - это решение системы линейных уравнений (8.1).
(8.6) - это решение системы линейных уравнений (8.2).
(8.7) - это решение системы линейных уравнений (8.3).
(8.8) - это решение системы линейных уравнений (8.4). □

**Теорема 8.2.** *Для любых значений параметров $a \neq 0$, $b \neq 0$, существует взаимно однозначное соответствие между координатами линейной функции алгебры $E(R, a, b)$ и её стандартными компонентами.*

*Доказательство.* Следствие теоремы 8.1. □

## 9. Регулярная функция

Хотя в алгебре кватернионов нет аналога уравнения Коши-Римана, в различных статьях изучаются различные множества функций, свойства которых похожи на свойства функций комплексного переменного. В [4, 5] определена регулярная функция, которая удовлетворяет уравнению

(9.1) $$\frac{\partial f}{\partial x^0} + i\frac{\partial f}{\partial x^1} + j\frac{\partial f}{\partial x^2} + k\frac{\partial f}{\partial x^3} = 0$$

**Теорема 9.1.** *Дифференциальное уравнение* (9.1) *эквивалентно системе дифференциальных уравнений*

(9.2) $$\begin{cases} \dfrac{\partial f^0}{\partial x^0} - \dfrac{\partial f^1}{\partial x^1} - \dfrac{\partial f^2}{\partial x^2} - \dfrac{\partial f^3}{\partial x^3} = 0 \\ \dfrac{\partial f^0}{\partial x^1} + \dfrac{\partial f^1}{\partial x^0} - \dfrac{\partial f^2}{\partial x^3} + \dfrac{\partial f^3}{\partial x^2} = 0 \\ \dfrac{\partial f^0}{\partial x^2} + \dfrac{\partial f^1}{\partial x^3} + \dfrac{\partial f^2}{\partial x^0} - \dfrac{\partial f^3}{\partial x^1} = 0 \\ \dfrac{\partial f^0}{\partial x^3} + \dfrac{\partial f^2}{\partial x^1} - \dfrac{\partial f^1}{\partial x^2} + \dfrac{\partial f^3}{\partial x^0} = 0 \end{cases}$$

*Доказательство.* Подставив

(9.3) $$\frac{\partial f}{\partial x^i} = \frac{\partial f^0}{\partial x^i} + \frac{\partial f^1}{\partial x^i}i + \frac{\partial f^2}{\partial x^i}j + \frac{\partial f^3}{\partial x^i}k$$

в уравнение (9.1), получим

(9.4) $$\begin{aligned}\frac{\partial f}{\partial x^0} + i\frac{\partial f}{\partial x^1} + j\frac{\partial f}{\partial x^2} + k\frac{\partial f}{\partial x^3} =& \frac{\partial f^0}{\partial x^0} - \frac{\partial f^1}{\partial x^1} - \frac{\partial f^2}{\partial x^2} - \frac{\partial f^3}{\partial x^3} \\ &+i(\frac{\partial f^0}{\partial x^1} + \frac{\partial f^1}{\partial x^0} - \frac{\partial f^2}{\partial x^3} + \frac{\partial f^3}{\partial x^2}) \\ &+j(\frac{\partial f^0}{\partial x^2} + \frac{\partial f^1}{\partial x^3} + \frac{\partial f^2}{\partial x^0} - \frac{\partial f^3}{\partial x^1}) \\ &+k(\frac{\partial f^0}{\partial x^3} + \frac{\partial f^2}{\partial x^1} - \frac{\partial f^1}{\partial x^2} + \frac{\partial f^3}{\partial x^0}) = 0\end{aligned}$$

Утверждение теоремы следует из равенства (9.4). □

В статье [1], следствие 3.1.2, p. 1000, показано, что если регулярная функция над алгеброй кватернионов имеет ограниченную норму, то это функция тождественно постоянна.

**Теорема 9.2.** *Компоненты производной Гато регулярной функции над алгеброй кватернионов удовлетворяют равенствам*

(9.5) $$\begin{aligned}\frac{\partial^{00}y}{\partial x} + \frac{\partial^{11}y}{\partial x} + \frac{\partial^{22}y}{\partial x} + \frac{\partial^{33}y}{\partial x} &= 0 \\ -\frac{\partial^{01}y}{\partial x} + \frac{\partial^{10}y}{\partial x} + \frac{\partial^{23}y}{\partial x} - \frac{\partial^{32}y}{\partial x} &= 0 \\ -\frac{\partial^{02}y}{\partial x} - \frac{\partial^{13}y}{\partial x} + \frac{\partial^{20}y}{\partial x} + \frac{\partial^{31}y}{\partial x} &= 0 \\ -\frac{\partial^{03}y}{\partial x} - \frac{\partial^{12}y}{\partial x} + \frac{\partial^{21}y}{\partial x} + \frac{\partial^{30}y}{\partial x} &= 0\end{aligned}$$

*Доказательство.* Подставив уравнения [2]-(7.4.1), [2]-(7.4.2), [2]-(7.4.3), [2]-(7.4.4) в уравнения (9.2), мы получим

$$\frac{\partial f^0}{\partial x^0} - \frac{\partial f^1}{\partial x^1} - \frac{\partial f^2}{\partial x^2} - \frac{\partial f^3}{\partial x^3}$$
$$= \frac{\partial^{00} y}{\partial x} - \frac{\partial^{11} y}{\partial x} - \frac{\partial^{22} y}{\partial x} - \frac{\partial^{33} y}{\partial x} - \left(\frac{\partial^{00} y}{\partial x} - \frac{\partial^{11} y}{\partial x} + \frac{\partial^{22} y}{\partial x} + \frac{\partial^{33} y}{\partial x}\right)$$
$$- \left(\frac{\partial^{00} y}{\partial x} + \frac{\partial^{11} y}{\partial x} - \frac{\partial^{22} y}{\partial x} + \frac{\partial^{33} y}{\partial x}\right) - \left(\frac{\partial^{00} y}{\partial x} + \frac{\partial^{11} y}{\partial x} + \frac{\partial^{22} y}{\partial x} - \frac{\partial^{33} y}{\partial x}\right)$$
$$= -2\frac{\partial^{00} y}{\partial x} - 2\frac{\partial^{11} y}{\partial x} - 2\frac{\partial^{22} y}{\partial x} - 2\frac{\partial^{33} y}{\partial x}$$
$$= 0$$

$$\frac{\partial f^0}{\partial x^1} + \frac{\partial f^1}{\partial x^0} - \frac{\partial f^2}{\partial x^3} + \frac{\partial f^3}{\partial x^2}$$
$$= -\frac{\partial^{01} y}{\partial x} - \frac{\partial^{10} y}{\partial x} + \frac{\partial^{23} y}{\partial x} - \frac{\partial^{32} y}{\partial x} + \left(\frac{\partial^{01} y}{\partial x} + \frac{\partial^{10} y}{\partial x} + \frac{\partial^{23} y}{\partial x} - \frac{\partial^{32} y}{\partial x}\right)$$
$$- \left(\frac{\partial^{01} y}{\partial x} - \frac{\partial^{10} y}{\partial x} - \frac{\partial^{23} y}{\partial x} - \frac{\partial^{32} y}{\partial x}\right) + \left(-\frac{\partial^{01} y}{\partial x} + \frac{\partial^{10} y}{\partial x} - \frac{\partial^{23} y}{\partial x} - \frac{\partial^{32} y}{\partial x}\right)$$
$$= -2\frac{\partial^{01} y}{\partial x} + 2\frac{\partial^{10} y}{\partial x} + 2\frac{\partial^{23} y}{\partial x} - 2\frac{\partial^{32} y}{\partial x}$$
$$= 0$$

$$\frac{\partial f^0}{\partial x^2} + \frac{\partial f^1}{\partial x^3} + \frac{\partial f^2}{\partial x^0} - \frac{\partial f^3}{\partial x^1}$$
$$= -\frac{\partial^{02} y}{\partial x} - \frac{\partial^{13} y}{\partial x} - \frac{\partial^{20} y}{\partial x} + \frac{\partial^{31} y}{\partial x} + \left(-\frac{\partial^{02} y}{\partial x} - \frac{\partial^{13} y}{\partial x} + \frac{\partial^{20} y}{\partial x} - \frac{\partial^{31} y}{\partial x}\right)$$
$$+ \left(\frac{\partial^{02} y}{\partial x} - \frac{\partial^{13} y}{\partial x} + \frac{\partial^{20} y}{\partial x} + \frac{\partial^{31} y}{\partial x}\right) - \left(\frac{\partial^{02} y}{\partial x} - \frac{\partial^{13} y}{\partial x} - \frac{\partial^{20} y}{\partial x} - \frac{\partial^{31} y}{\partial x}\right)$$
$$= -2\frac{\partial^{02} y}{\partial x} - 2\frac{\partial^{13} y}{\partial x} + 2\frac{\partial^{20} y}{\partial x} + 2\frac{\partial^{31} y}{\partial x}$$
$$= 0$$

$$\frac{\partial f^0}{\partial x^3} + \frac{\partial f^2}{\partial x^1} - \frac{\partial f^1}{\partial x^2} + \frac{\partial f^3}{\partial x^0}$$
$$= -\frac{\partial^{03} y}{\partial x} + \frac{\partial^{12} y}{\partial x} - \frac{\partial^{21} y}{\partial x} - \frac{\partial^{30} y}{\partial x} - \left(-\frac{\partial^{03} y}{\partial x} - \frac{\partial^{12} y}{\partial x} - \frac{\partial^{21} y}{\partial x} + \frac{\partial^{30} y}{\partial x}\right)$$
$$- \left(\frac{\partial^{03} y}{\partial x} - \frac{\partial^{12} y}{\partial x} - \frac{\partial^{21} y}{\partial x} - \frac{\partial^{30} y}{\partial x}\right) - \left(\frac{\partial^{03} y}{\partial x} + \frac{\partial^{12} y}{\partial x} - \frac{\partial^{21} y}{\partial x} + \frac{\partial^{30} y}{\partial x}\right)$$
$$= -2\frac{\partial^{03} y}{\partial x} - 2\frac{\partial^{12} y}{\partial x} + 2\frac{\partial^{21} y}{\partial x} + 2\frac{\partial^{30} y}{\partial x}$$
$$= 0$$

□

**Теорема 9.3.** *Дифференциал Гато регулярной функции над алгеброй кватернионов имеет вид*

(9.6)
$$-\left(\frac{\partial^{11} y}{\partial x} + \frac{\partial^{22} y}{\partial x} + \frac{\partial^{33} y}{\partial x}\right) dx + \frac{\partial^{11} y}{\partial x} i dx i + \frac{\partial^{22} y}{\partial x} j dx j + \frac{\partial^{33} y}{\partial x} k dx k$$
$$+ \left(\frac{\partial^{10} y}{\partial x} + \frac{\partial^{23} y}{\partial x} - \frac{\partial^{32} y}{\partial x}\right) dx i + \frac{\partial^{10} y}{\partial x} i dx + \frac{\partial^{23} y}{\partial x} j dx k - \frac{\partial^{32} y}{\partial x} k dx i$$
$$+ \left(-\frac{\partial^{13} y}{\partial x} + \frac{\partial^{20} y}{\partial x} + \frac{\partial^{31} y}{\partial x}\right) dx j - \frac{\partial^{13} y}{\partial x} i dx k + \frac{\partial^{20} y}{\partial x} j dx + \frac{\partial^{31} y}{\partial x} k dx i$$
$$+ \left(-\frac{\partial^{12} y}{\partial x} + \frac{\partial^{21} y}{\partial x} + \frac{\partial^{30} y}{\partial x}\right) dx k - \frac{\partial^{12} y}{\partial x} i dx j + \frac{\partial^{21} y}{\partial x} j dx i + \frac{\partial^{30} y}{\partial x} k dx$$

*Доказательство.* Утверждение теоремы является следствием теоремы [2]-5.2.8. □

Равенство (9.1) эквивалентно равенству

(9.7)
$$\begin{pmatrix} 1 & i & j & k \end{pmatrix} *^* \begin{pmatrix} \dfrac{\partial f^0}{\partial x^0} & \dfrac{\partial f^0}{\partial x^1} & \dfrac{\partial f^0}{\partial x^2} & \dfrac{\partial f^0}{\partial x^3} \\ \dfrac{\partial f^1}{\partial x^0} & \dfrac{\partial f^1}{\partial x^1} & \dfrac{\partial f^1}{\partial x^2} & \dfrac{\partial f^1}{\partial x^3} \\ \dfrac{\partial f^2}{\partial x^0} & \dfrac{\partial f^2}{\partial x^1} & \dfrac{\partial f^2}{\partial x^2} & \dfrac{\partial f^2}{\partial x^3} \\ \dfrac{\partial f^3}{\partial x^0} & \dfrac{\partial f^3}{\partial x^1} & \dfrac{\partial f^3}{\partial x^2} & \dfrac{\partial f^3}{\partial x^3} \end{pmatrix} *^* \begin{pmatrix} 1 \\ i \\ j \\ k \end{pmatrix} = 0$$

## 10. Вместо эпилога

Поле комплексных чисел и алгебра кватернионов имеют и общие свойства, и некоторые различия. Эти различия делают более трудной задачу найти в алгебре кватернионов закономерности, подобные знакомым нам в поле комплексных чисел. Поэтому понимание этих различий очень важно.

Одно из направлений исследования - это найти аналог уравнения Коши-Римана в алгебре кватернионов. В статье я проанализировал некоторые исследования в этой области.

Согласно теореме [2]-7.1.1 линейное отображение имеет матрицу

$$\begin{pmatrix} a_0 & -a_1 \\ a_1 & a_0 \end{pmatrix}$$

Это отображение соответствует умножению на число $a = a_0 + a_1 i$. Утверждение следует из равенств

$$(a_0 + a_1 i)(x_0 + x_1 i) = a_0 x_0 - a_1 x_1 + (a_0 x_1 + a_1 x_0) i$$
$$\begin{pmatrix} a_0 & -a_1 \\ a_1 & a_0 \end{pmatrix} \begin{pmatrix} x_0 \\ x_1 \end{pmatrix} = \begin{pmatrix} a_0 x_0 - a_1 x_1 \\ a_1 x_0 + a_0 x_1 \end{pmatrix}$$

Я решил рассмотреть аналогичный класс функций кватернионов. Линейное отображение алгебры кватернионов

$$x \to ax$$

имеет матрицу

$$
(10.1) \quad \begin{pmatrix} a^0 & -a^1 & -a^2 & -a^3 \\ a^1 & a^0 & -a^3 & a^2 \\ a^2 & a^3 & a^0 & -a^1 \\ a^3 & -a^2 & a^1 & a^0 \end{pmatrix}
$$

Интересно рассмотреть класс функций кватерниона, которые имеют производную похожей структуры. Однако структура матрицы (10.1) мне представляется несколько жёсткой для производной функции кватернионов и я несколько ослаблю требование. Я буду предполагать, что производная удовлетворяет следующим равенствам

$$
(10.2) \quad \frac{\partial y^0}{\partial x^0} = \frac{\partial y^1}{\partial x^1} = \frac{\partial y^2}{\partial x^2} = \frac{\partial y^3}{\partial x^3}
$$

$$
(10.3) \quad \frac{\partial y^i}{\partial x^j} = -\frac{\partial y^j}{\partial x^i} \quad i \neq j
$$

Нетрудно убедиться, что производная функции вида

$$y = ax \quad y = xa$$

удовлетворяет равенствам (10.2), (10.3).

Рассмотрим функцию

$$y = x^2$$

Непосредственное вычисление даёт

$$
(10.4) \quad \begin{aligned} y^0 &= (x^0)^2 - (x^1)^2 - (x^2)^2 - (x^3)^2 \\ y^1 &= 2x^0 x^1 \\ y^2 &= 2x^0 x^2 \\ y^3 &= 2x^0 x^3 \end{aligned}
$$

Производная отображения (10.4) имеет матрицу

$$
(10.5) \quad \begin{pmatrix} 2x^0 & -2x^1 & -2x^2 & -2x^3 \\ 2x^1 & 2x^0 & 0 & 0 \\ 2x^2 & 0 & 2x^0 & 0 \\ 2x^3 & 0 & 0 & 2x^0 \end{pmatrix}
$$

Следовательно, матрица (10.5) удовлетворяет равенствам (10.2), (10.3).

**Теорема 10.1.** *Производная сопряжения не удовлетворяет равенствам* (10.2), (10.3).

*Доказательство.* Согласно доказательству теоремы [2]-7.3.5, производная сопряжения не удовлетворяет равенству (10.2). □

Однако множество функций, производная которых удовлетворяет равенствам (10.2), (10.3) не достаточно велико.

**Теорема 10.2.** *Отображение*

$$(10.6) \quad y = x^3$$

*не удовлетворяет равенству* (10.2).

*Доказательство.* Мы можем отображение (10.6) представить в виде

$$y = x\,x^2$$

Координаты этого отображения имеют вид

$$
\begin{aligned}
(10.7)\quad y^0 &= x^0((x^0)^2 - (x^1)^2 - (x^2)^2 - (x^3)^2) \\
&\quad - x^1 2x^0 x^1 - x^2 2x^0 x^2 - x^3 2x^0 x^3 \\
&= (x^0)^3 - x^0(x^1)^2 - x^0(x^2)^2 - x^0(x^3)^2 \\
&\quad - 2x^0(x^1)^2 - 2x^0(x^2)^2 - 2x^0(x^3)^2 \\
&= (x^0)^3 - 3x^0(x^1)^2 - 3x^0(x^2)^2 - 3x^0(x^3)^2
\end{aligned}
$$

$$
\begin{aligned}
(10.8)\quad y^1 &= x^1((x^0)^2 - (x^1)^2 - (x^2)^2 - (x^3)^2) \\
&\quad + x^0 2x^0 x^1 - x^3 2x^0 x^2 + x^2 2x^0 x^3 \\
&= x^1(x^0)^2 - (x^1)^3 - x^1(x^2)^2 - x^1(x^3)^2 \\
&\quad + 2(x^0)^2 x^1 - 2x^0 x^2 x^3 + 2x^0 x^2 x^3 \\
&= 3x^1(x^0)^2 - (x^1)^3 - x^1(x^2)^2 - x^1(x^3)^2
\end{aligned}
$$

Из равенства (10.7), следует

$$(10.9)\qquad \frac{\partial y^0}{\partial x^0} = 3(x^0)^2 - 3(x^1)^2 - 3(x^2)^2 - 3(x^3)^2$$

Из равенства (10.8), следует

$$(10.10)\qquad \frac{\partial y^1}{\partial x^1} = 3(x^0)^2 - 3(x^1)^2 - (x^2)^2 - (x^3)^2$$

Утверждение теоремы следует из равенств (10.9), (10.10). $\square$

Без сомнения, это только начало исследования, и немало вопросов должно быть отвечено.

## 11. Список литературы